\magnification=\magstep1\baselineskip 12pt

\font\Bbb=msbm10

\centerline{{\bf On the Brauer groups of quasilocal fields and}}
\smallskip \centerline{\bf the norm groups of their finite Galois
extensions\footnote{$^{\ast }$}{\rm 2000 Mathematics Subject
Classification: Primary 12F10, 16K50; Secondary 12J25, 14F22}}
\vskip0.6truecm

\centerline{I.D. Chipchakov\footnote{$^{\ast}$}{Partially
supported by Grant MI-1503/2005 of the Bulgarian National Science Fund.}}
\par
\medskip
\centerline{Institute of Mathematics and Informatics, Bulgarian
Academy of Sciences}
\par
\centerline{Acad. G. Bonchev Str., bl. 8, 1113 Sofia, Bulgaria;
chipchak@math.bas.bg}
\par
\vskip0.9truecm Abstract. This paper shows that divisible abelian
torsion groups are realizable as Brauer groups of quasilocal fields.
It describes the isomorphism classes of Brauer groups of primarily 
quasilocal fields and solves the analogous problem concerning the 
reduced components of the Brauer groups of two basic types of 
Henselian valued absolutely stable fields. For a quasilocal field 
$E$ and a finite separable extension $R/E$, we find two sufficient 
conditions for validity of the norm group equality $N(R/E) = N(R 
_{0}/E)$, where $R _{0}$ is the maximal abelian extension of $E$ in 
$R$. This is used for deriving information on the arising specific 
relations between Galois groups and norm groups of finite Galois 
extensions of $E$.
\par
\medskip
Key words: Quasilocal field; Brauer group; Character group;
Corestriction; Galois extension; Norm group; Abelian closed class;
Transfer; Brauer-Severi variety.
\par
\vskip0.9truecm \centerline{{\bf 1. Introduction and statements of
the main results}}
\par
\medskip
This paper is devoted to the study of norm groups and Brauer groups
of the fields pointed out in the title, i.e. of fields whose finite
extensions are primarily quasilocal (abbr., PQL). Our main result
describes, up-to an isomorphism, the abelian groups that can be
realized as Brauer groups of several basic types of PQL-fields (see
Theorem 1.2, Propositions 2.3 (ii), 3.4 and Section 6). For a
quasilocal field $E$ and a finite separable extension $R/E$, it 
gives two sufficient conditions that the norm group $N(R/E)$ 
coincides with $N(R/E) _{\rm Ab}$, the norm group of the maximal
abelian extension of $E$ in $R$ (see Theorem 1.1). When the field 
$E$ is nonreal, this allows us to clarify essential algebraic and 
topological aspects of the behaviour of norm groups of finite Galois 
extensions of $E$.
\par
\medskip
The basic notions needed to present this research are the same as
those in [7]; the reader is referred to [18; 21; 27; 31 and 15], for
any missing definitions concerning simple algebras, Brauer groups,
field extensions, Galois cohomology and abelian groups. Simple
algebras are supposed to be associative with a unit and
finite-dimensional over their centres, and Galois groups are viewed
as profinite with respect to the Krull topology. For a central
simple algebra $A$ over a field $E$, we write $[A]$ for the
similarity class of $A$ in the Brauer group Br$(E)$. As usual, $E
^{\ast }$ denotes the multiplicative group of $E$, $E _{\rm sep}$ a
separable closure of $E$, $G _{E} = G(E _{\rm sep}/E)$ is the
absolute Galois group of $E$, $C _{E}$ stands for the character
group of $G _{E}$, and $d(E)$ is the class of central division
$E$-algebras. For an arbitrary field extension $\Lambda /E$,
Br$(\Lambda /E)$ is its relative Brauer group, $\rho _{E/\Lambda }$
is the scalar extension map, Im$(E/\Lambda )$ is the image of $\rho
_{E/\Lambda }$, and $I(\Lambda /E)$ is the set of intermediate
fields of $\Lambda /E$. We write Gal$(E)$ for the set of finite
Galois extensions of $E$ in $E _{\rm sep}$, and put $\Omega (E) =
\{M \in {\rm Gal}(E)\colon \ G(M/E) \in {\rm Ab}\}$. Throughout,
$\overline P$ is the set of prime numbers and $\Pi (E)$ consists of
those $p \in \overline P$, for which $G _{E}$ is of nonzero
cohomological $p$-dimension cd$_{p} (G _{E})$. For each $p \in
\overline P$, $E (p)$ is the maximal $p$-extension of $E$ in $E
_{\rm sep}$, and $\Omega _{p} (E) = \{Y \in \Omega (E)\colon \ Y
\subseteq E (p)\}$. When $\Psi $ is a nonempty formation of finite
groups in the sense of [34], $E _{\Psi }$ denotes the compositum of
all fields $M \in {\rm Gal}(E)$ with $G(M/E) \in \Psi $; in view of
Galois theory and the choice of $\Psi $, $E _{\Psi }$ is the union
of these $M$. For a finite extension $R$ of $E$ in $E _{\rm sep}$,
we put $N(R/E) _{\Psi } = N(E _{\Psi } \cap R/E)$. The formations
of abelian, metabelian, nilpotent, solvable, and of all finite
groups are denoted by ${\rm Ab}$, ${\rm Met}$, ${\rm Nil}$, ${\rm
Sol}$ and ${\rm Fin}$, respectively. A class $\chi \subseteq {\rm
Fin}$ is called abelian closed, if it is nonempty and closed with
respect to taking subgroups, homomorphic images, finite direct
products, and group extensions with abelian kernels (a series of
typical examples of such classes is given in Remark 6.1). We say
that $E$ is formally real, if $-1$ is not presentable over $E$ as a
finite sum of squares; $E$ is called nonreal, otherwise. The field
$E$ is said to be PQL, if every cyclic extension $F$ of $E$ embeds
as an $E$-subalgebra in each $D \in d(E)$ of Schur index ind$(D)$
divisible by the degree $[F\colon E]$. We say that $E$ is strictly
PQL, if it is PQL and the $p$-component Br$(E) _{p}$ of Br$(E)$ is
nontrivial, when $p$ runs through the set $P(E)$ of those elements
of $\overline P$, for which $E (p) \neq E$.
\par
\medskip
Singled out in the process of characterizing basic types of stable
fields with Henselian valuations (see [7] and the references there),
PQL-fields $E$ deserve interest in their own right because of the
arising close relations between the fields in $\Omega (E)$, their
norm groups and central simple $E$-algebras. Firstly, it should be
pointed out that strictly PQL-fields admit one-dimensional local
class field theory (abbr. LCFT, see (2.1)) and the converse holds in
all presently known cases (cf. [8, Theorem 1.1, Remark 4.4 and Sect.
3]). Note also that the field $E$ is strictly quasilocal (SQL), i.e.
its finite extensions are strictly PQL, if and only if they admit
LCFT [8, Proposition 3.6]. Secondly, this research is motivated by the
dependence of some of these relations on the structure of Br$(E)$
(see [8, Theorems 1.2 and 1.3]), and therefore, by the problem of
describing the isomorphism classes of Brauer groups over the main
kinds of PQL-fields. It is worth adding that the quasilocal property
singles out one of the basic classes of absolutely stable fields (in
the sense of Brussel, see [7, I, Proposition 2.3]), and the
structure of Br$(F)$, for an arbitrary absolutely stable field $F$,
is of interest for the theory of central simple algebras in general
(see [27, Sects. 14.4 and 19.6]). The choice of our main topic is
determined by the fact that the groups $N(M/E)\colon \ M \in {\rm
Gal}(E) \setminus \Omega (E)$, reflect more aspects of the specific
nature of $E$ than merely the influence of the PQL-property. Our
starting point are the following analogues to the norm limitation
theorem about local fields (see [17, Ch. 6, Theorem 8]):
\par
\medskip
(1.1) (i) $N(R/E) = N(R/E) _{\rm Ab}$, provided that $R$ is a finite
separable extension of a field $E$ with LCFT in the sense of
Neukirch-Perlis [26], i.e. if the triple $(G _{E}, \{G(E _{\rm
sep}/F),$ $F \in {\rm Fe}(E)\}, E _{\rm sep} ^{\ast })$ is an
Artin-Tate class formation (cf. [2, Ch. XIV]), where Fe$(E)$ is the
set of finite extensions of $E$ in $E _{\rm sep}$;
\par
(ii) $N(R/E) = N(R/E) _{\rm Ab}$, if $E$ is PQL and $R \subseteq E
_{\rm Nil}$ [5].
\par
\medskip
It is known (cf. [26]) that a field $E$ admits LCFT in the sense of
Neukirch-Perlis if and only if it is SQL, Br$(E)$ embeds in the
quotient group $\hbox{\Bbb Q}/\hbox{\Bbb Z}$ of the additive group
of rational numbers by the subgroup of integers, and $\rho
_{E/F}\colon {\rm Br}(E) \to {\rm Br}(F)$ is surjective, for every
finite extension $F$ of $E$ in $E _{\rm sep}$. This holds when $E$
has a Henselian discrete valuation with a quasifinite residue field
$\widehat E$ (see, e.g., [41]). The basis for the present discussion
is also formed by the characterization of the PQL-property in the
class of algebraic extensions of global fields, which yields the
following (see [5, Sects. 1 and 2] and the references there):
\par
\medskip
(1.2) (i) For each $G \in {\rm Fin} \setminus {\rm Nil}$, there
exist algebraic extensions $E(G)$ and $M(G)$ of the field
$\hbox{\Bbb Q}$ of rational numbers, such that $E(G)$ is strictly
PQL, $M(G) \in {\rm Gal}(E(G))$, $G(M(G)/E(G)) \cong G$ and
$N(M(G)/E(G)) \neq N(M(G)/E(G)) _{\rm Ab}$;
\par
(ii) If $E$ is an algebraic PQL-extension of a global field $E
_{0}$, then Br$(E)$ embeds in $\hbox{\Bbb Q}/\hbox{\Bbb Z}$.
Moreover, if $R/E$ is a finite extension, then $N(R/E) = N(\Sigma
/E)$, for some $\Sigma \in \Omega (E)$; when $E$ is strictly PQL,
$\Sigma $ is uniquely determined by $R/E$.
\par
\medskip
The purpose of this paper is to present two main results which shed
an additional light on (1.1) and (1.2), and solve the above-noted
problem for Brauer groups of nonreal PQL-fields. The first result is
stated as follows:
\par
\medskip
{\bf Theorem 1.1.} {\it Let $E$ be a quasilocal field and $R$ a
finite extension of $E$ in $E _{\rm sep}$. Then $N(R/E) = N(R/E)
_{\rm Ab}$ in the following two cases:}
\par
(i) {\it The map $\rho _{E/M}$ is surjective, for some $M \in {\rm
Gal}(E)$ including $R$;}
\par
(ii) {\it There exists a field $\Phi (R) \in \Omega (E)$, such that
$N(\Phi (R)/E) \subseteq N(R/E)$.}
\par
\medskip
Theorem 1.1 is deduced in Section 3 from its $p$-primary analogue
stated as Theorem 3.1. This analogue enables us to generalize
Theorem 1.1 (i) by proving at the end of Section 3 that $N(R/E) =
N(R/E) _{\rm Ab}$, provided that $E$ is quasilocal and $R \in
I(M/E)$, for some $M \in {\rm Gal}(E)$ with $\rho _{E/L}$
surjective, where $L$ is the fixed field of the Fitting subgroup of
$G(M/E)$. In this setting, it may occur that $\rho _{E/\Phi }$ is
not surjective, for any $\Phi \in {\rm Gal}(E)$ including $R$ (see
the comment preceding Proposition 6.3). Theorem 3.1 has been used in
[6] for describing the norm groups of finite separable extensions of
SQL-fields with Henselian discrete valuations. Like the description
of the norm groups of formally real quasilocal fields, obtained in
[9], this yields a generally nonabelian LCFT. Our second main
result, combined with [7, I, Theorem 3.1 (ii) and Lemma 3.5], shows
that an abelian torsion group $T$ is isomorphic to Br$(E)$, for some
nonreal PQL-field $E$, if and only if $T$ is divisible (for the
formally real case, see Propositions 6.4 (i) and 3.4); this
specifies observations made at the end of [37, Sect. 3]. When $T$ is
divisible, it states that $E$ can be found among quasilocal fields
so as to solve one of the main inverse problems related to (1.1) and
(1.2).
\par
\medskip
{\bf Theorem 1.2.} {\it Let $E _{0}$ be a field, $T$ a divisible
abelian torsion group, $T _{0}$ a subgroup of {\rm Br}$(E _{0})$
embeddable in $T$, and let $\chi $ and $\chi ^{\prime }$ be
subclasses of ${\rm Fin}$, such that ${\rm Nil} \subseteq \chi
\subseteq \chi ^{\prime }$. Assume also that any class $\chi $,
$\chi ^{\prime }$ is abelian closed unless it equals ${\rm Nil}$.
Then there exists a quasilocal and nonreal extension $E = E(T)$ of
$E _{0}$ with the following properties:}
\par
(i) Br$(E) \cong T${\it , $E _{0}$ is separably closed in $E$, $\rho
_{E/E _{0}}$ maps $T _{0}$ injectively into {\rm Br}$(E)$, and each
$G \in {\rm Fin}$ is realizable as a Galois group over $E$;}
\par
(ii) {\it For each finite extension $R$ of $E$ in $E _{\chi }$,
$N(R/E) = N(R/E) _{\rm Ab}$; moreover, if $\chi \neq {\rm Nil}$,
then $\rho _{E/R}$ is surjective;}
\par
(iii) $N(M/E) \neq N(M/E) _{\rm Ab}${\it , for every $M \in {\rm
Gal}(E)$ with $G(M/E) \not\in \chi ^{\prime }$;}
\par
(iv) {\it If $G \in \chi ^{\prime } \setminus \chi $, then there
are $M _{1}, M _{2} \in {\rm Gal}(E)$, such that $G(M _{j}/E)$
$\cong G$, $j =  1, 2$, $N(M _{1}/E) = N(M _{1}/E) _{\rm Ab}$ and
$N(M _{2}/E) \neq N(M _{2}/E) _{\rm Ab}$.}
\par
\medskip
The assertions of Theorem 1.2 (i)-(ii) in the case of $T _{0} = {\rm
Br}(E _{0})$, combined with [7, I, Corollary 8.5] and the behaviour
of Schur indices under scalar extensions of finite degrees (cf. [27,
Sect. 13.4]), imply [37, Theorems 3.7$\div $3.9]. Since
$n$-dimensional $F$-algebras embed in the matrix $F$-algebra $M _{n}
(F)$, for any field $F$ and $n \in \hbox{\Bbb N}$, these assertions
and well-known properties of tensor products (see [27, Sects. 9.3
and 9.4 ]) also enable one to deduce [37, Theorem 3.10] from the
Skolem-Noether theorem (as in the proof of the double centralizer
theorem, for example, in [27, Sect. 12.7]). Thus the noted part of
Theorem 1.2 simplifies the proofs of [37, Theorems 2.6 and 2.8].
When $\chi = {\rm Fin}$, it admits a Galois cohomological
interpretation (see Remark 5.4 and [7, I, Theorem 8.1]) and yields
the following result:
\par
\medskip
(1.3) For each divisible abelian torsion group $T$, there is a
quasilocal field $E$, such that Br$(E) \cong T$, all $G \in {\rm
Fin}$ can be realized as Galois groups over $E$, and $N(R/E) =
N(R/E) _{\rm Ab}$, for every finite extension $R$ of $E$ in $E _{\rm
sep}$.
\par
\medskip
When $\chi = {\rm Nil}$ or $\chi $ is abelian closed with ${\rm Nil}
\subset \chi \neq {\rm Fin}$, the conclusions of Theorem 1.2 in the
special cases of $\chi ^{\prime } = \chi $ and  $\chi ^{\prime } =
{\rm Fin}$ amount essentially to the following:
\par
\medskip
(1.4) For each divisible abelian torsion group $T$, there exist
quasilocal fields $E _{1}$ and $E _{2}$ with Br$(E _{i}) \cong T$,
$i = 1, 2$, and such that:
\par
(i) All $G \in {\rm Fin}$ are realizable as Galois groups over $E
_{1}$ and $E _{2}$, and whenever $M _{1} \in {\rm Gal}(E _{1})$,
$N(M _{1}/E _{1}) = N(M _{1}/E _{1}) _{\rm Ab}$ if and only if $G(M
_{1}/E _{1}) \in \chi $.
\par
(ii) $N(M _{2}/E _{2}) = N(M _{2}/E _{2}) _{\rm Ab}$, provided that
$M _{2} \in {\rm Gal}(E _{2})$ and $G(M _{2}/E _{2}) \in \chi $. For
each $G \in {\rm Fin} \setminus \chi $, Gal$(E _{2})$ contains
elements $M(G) _{1}$ and $M(G) _{2}$ with
$G(M(G) _{j}/E _{2}) \cong G$, $j = 1, 2$, $N(M(G) _{1}/E _{2}) =
N(M(G) _{1}/E _{2}) _{\rm Ab}$ and $N(M(G) _{2}/E _{2})$ $\neq
N(M(G) _{2}/E _{2}) _{\rm Ab}$.
\par
\medskip
Theorem 1.1 and statements (1.1) (ii), (1.2) (i)-(ii), (1.3) and
(1.4) mark the limit behaviour of norm groups of finite Galois
extensions of PQL-fields. By [8, (2.3)], the fields singled out by
(1.3) and (1.4) have no Henselian valuations with indivisible value
groups. Note also that if $(F, v)$ is a Henselian discrete valued
SQL-field, then $G _{F}$ is prosolvable of special type (see
Corollary 6.7 and [8, (2.1) and the comments to (2.4) (ii)]). These 
facts and the topological interpretation of Theorem 1.2 in Section 
6 allow one to appreciate from an algebraic point of view the 
Neukirch-Perlis generalization of LCFT, and without artificial 
limitations, to incorporate it in the study of quasilocal fields and 
other areas.
\par
\medskip
Here is an overview of the paper: Section 2 includes preliminaries
needed in the sequel, such as statements of frequently used projection
formulae relating the corestrictions of Brauer and character groups
of an arbitrary finite separable extension. The proofs of these
formulae (and of Proposition 2.8) given in [10] as well as of
Propositions 4.1 and 6.8 show that Theorem 1.2 (i)-(ii), applied $T
_{0} = {\rm Br}(E _{0})$ and $\chi = {\rm Fin}$, provides useful
tools for the study of various aspects of Brauer group theory on a
unified basis. Theorems 1.1 and 1.2 are proved in Sections 3 and 5,
respectively. The technical preparation for the proof of Theorem 1.2
is made in Section 4 and its results seem to be of independent
interest. As an application of Theorem 1.2, we describe in Section 6
the isomorphism classes of Brauer groups of formally real PQL and of
strictly PQL-fields, and do the same for the reduced parts of the
Brauer groups of two basic types of Henselian valued absolutely
stable fields.
\par
\vskip0.6truecm \centerline {\bf 2. Preliminaries on norm groups,
$p$-quasilocal fields}
\smallskip
\centerline{\bf and corestrictions of Brauer and character groups}
\par
\medskip
(2.1) Let $E$ be a field and Nr$(E)$ the set of norm groups of its
finite extensions in $E _{\rm sep}$. We say that $E$ admits
LCFT, if the mapping $\pi \colon \ \Omega (E) \to {\rm Nr}(E)$, by
the rule $\pi (F) = N(F/E)\colon \ F \in \Omega (E)$, is injective,
and whenever $M _{1}, M _{2} \in \Omega (E)$, $N(M _{1}M _{2}/E) =
N(M _{1}/E) \cap N(M _{2}/E)$ and $N(M _{1} \cap M _{2}/E) = N(M
_{1}/E)N(M _{2}/E)$ (as usual, $M _{1}M _{2}$ is the compositum of
$M _{1}$ and $M _{2}$). We call $E$ a field with local $p$-class
field theory (local $p$-CFT), for some $p \in \overline P$, if the
restriction of $\pi $ on the set $\Omega _{p} (E)$ has the same
properties.
\par
\medskip
The following lemma (proved, e.g., in [5]) implies that a field $E$
admits LCFT if and only if it admits local $p$-CFT, for every $p \in
P(E)$. When $E$ is of this kind, [7, I, Lemma 4.2 (ii)] shows that
Br$(E) _{p} \neq \{0\}$, $p \in P(E)$.
\par
\medskip
{\bf Lemma 2.1.} {\it Let $E$, $R$ and $M$ be fields, such that $R
\in I(M/E)$, $R \neq E$, $M \in {\rm Gal}(E)$ and $G(M/E) \in {\rm
Nil}$. Let $P(R/E)$ be the set of prime divisors of $[R\colon E]$,
and $R _{p} = R \cap E (p)$, for each $p \in P(R/E)$. Then $R$
equals the compositum of the fields $R _{p}\colon \ p \in P(R/E)$,
$[R\colon E] = \prod _{p \in P(R/E)} [R _{p}\colon E]$, $N(R/E) =
\cap _{p \in P(R/E)} N(R _{p}/E)$ and the quotient group $E ^{\ast
}/N(R/E)$ is isomorphic to the direct group product $\prod _{p \in
P(R/E)} E ^{\ast }/N(R _{p}/E)$.}
\par
\medskip
Henceforth, Syl$_{p} (M/E)$ denotes the set of Sylow $p$-subgroups
of $G(M/E)$, for any $M \in {\rm Gal}(E)$, $p \in \overline P$. For
the proof of the following lemma, we refer the reader to [7, II].
\par
\medskip
{\bf Lemma 2.2.} {\it Let $E$ and $M$ be fields, $M \in {\rm
Gal}(E)$ and $P(M/E) = \{p \in \overline P\colon \ p \vert [M\colon
E]\}$. Then $N(M/E) \subseteq N(M/F)$, for every $F \in I(M/E)$.
Moreover, if $E _{p}$ is the fixed field of a group $G _{p} \in
{\rm Syl}_{p} (M/E)$, then $N(M/E) = \cap _{p \in P(M/E)} N(M/E
_{p})$.}
\par
\medskip
The main results of [7, I] used in the present paper (supplemented
by a well-known result on orderings in Pythagorean fields), can be
stated as follows:
\par
\medskip
{\bf Proposition 2.3.} {\it Assume that $E$ is a $p$-quasilocal
field, for some $p \in P(E)$, $R$ is a finite extension of $E$ in $E
(p)$, and $D \in d(E)$ is an algebra of $p$-primary index. Then:}
\par
(i) $R$ {\it is a $p$-quasilocal field and {\rm ind}$(D) = {\rm
exp}(D)$.}
\par
(ii) Br$(R) _{p}$ {\it is a divisible group unless $p = 2$, $R = E$
and $E$ is formally real. In the formally real case, $E (2) =
E(\sqrt{-1})$, {\rm Br}$(E) _{2}$ is of order $2$ and {\rm
Br}$(E(\sqrt{-1})) _{2} = \{0\}$; this occurs if and only if $E$ is
Pythagorean with a unique ordering.}
\par
(iii) $\rho _{E/R}$ {\it maps {\rm Br}$(E) _{p}$ surjectively on
{\rm Br}$(R) _{p}$, i.e. {\rm Br}$(R) _{p} \subseteq {\rm Im}(E/R)$.}
\par
(iv) {\it $R$ embeds in $D$ as an $E$-subalgebra if and only if
$[R\colon E] \vert {\rm ind}(D)$.}
\par
\medskip
The following lemma provides an easy method of constructing
$p$-quasilocal fields. Before stating it, recall that a field
extension $F/F _{0}$ is said to be regular, if $F _{0}$ is separably
closed in $F$ and $I(F/F _{0})$ contains an element $F _{0} ^{\prime
}$, such that $F _{0} ^{\prime }/F _{0}$ is rational (i.e. purely
transcendental) and $F/F _{0} ^{\prime }$ is separable. It is known
that the tensor product $\otimes _{F _{0}} \Lambda _{i}$, $i \in I$,
of such extensions is a domain with a fraction field $F(I)$ regular
over $F _{0}$. We call $\Lambda (I)$ a tensor compositum of the
fields $\Lambda _{i}$ over $F _{0}$, and write $\Lambda (I) =
\widetilde \otimes _{F _{0}} \Lambda _{i}$. Recall further that the
class Reg$(F _{0})$ of regular extensions of $F _{0}$ contains the
function fields of the $F _{0}$-varieties (i.e. algebraic varieties
defined over $F _{0}$ and irreducible over $F _{0,{\rm sep}}$)
considered in this paper. In what follows, we shall use without an
explicit reference the well-known facts that $F _{1} \otimes _{F
_{0}} \Lambda \in {\rm Reg}(F _{1})$ whenever $\Lambda \in {\rm
Reg}(F _{0})$ and $F _{1}$ is a finite extension of $F _{0}$ in $F
_{0,{\rm sep}}$, the image of Reg$(F _{1})$ under the transfer map
Tr$_{F _{1}/F _{0}}$ (over $F _{1}/F _{0}$) is included in Reg$(F
_{0})$, and the compositions $\widetilde \otimes _{F _{0}} \circ
{\rm Tr}_{F _{1}/F _{0}}$ and Tr$_{F _{1}/F _{0}} \circ \widetilde
\otimes _{F _{1}}$ coincide. These are easily obtained from Galois
theory and the definition of Tr$_{F _{1}/F _{0}}$ (see the beginning
of [30, Sect. 3]).
\par
\medskip
{\bf Lemma 2.4.} {\it Let $F _{0}$ be a field and $p \in \overline
P$. Then there exists a field extension $F/F _{0}$, such that $F$ is
$p$-quasilocal, $F _{0}$ is algebraically closed in $F$ and {\rm
Br}$(F/F _{0}) = \{0\}$.}
\par
\medskip
{\it Proof.} Using [14, Theorem 1], one constructs $F$ as a
union $\cup _{i=1} ^{\infty } B _{i} = \cup _{i=1} ^{\infty } B
_{i} ^{\prime }$ of fields defined inductively as follows:
\par
\medskip
(2.2) (i) $B _{0}$ is a rational function field in one
indeterminate over $F _{0}$ and $B _{0} ^{\prime } = B _{0}$
(since $p \in P(B _{0})$, this implies that Br$(B _{0} ^{\prime })
_{p} \neq \{0\}$).
\par
(ii) For each $i \in \hbox{\Bbb N}$, $B _{i}/B _{(i-1)} ^{\prime
}$ is a rational extension with a transcendence basis (abbr,
tr-basis) $\{X _{(r _{i},c _{i})}\}$ indexed by the Cartesian
product $R _{p,i} \times C _{p,i}$, where $R _{p,i}$ is a system of representatives of the isomorphism classes of algebras in $d(B
_{(i-1)} ^{\prime })$ of index $p$, and $C _{p,i}$ is the set of
extensions of $B _{(i-1)} ^{\prime }$ in $B _{(i-1)} ^{\prime } (p)$
of degree $p$.
\par
(iii) For any $i \in \hbox{\Bbb N}$, let $F(r _{i}, c _{i})$ be the 
function field of the Brauer-Severi $B _{i}$-variety canonically
associated with the central simple $B _{i}$-algebra $(r _{i} \otimes
$ $_{B _{(i-1)}'} B _{i}) \otimes ((c _{i} \otimes _{B _{(i-1)}'} B
_{i})/B _{i}, \tilde c _{i}, X _{(r _{i},c _{i})} ^{-1})$, for
each $(r _{i}, c _{i})$ $\in R _{p,i} \times C _{p,i}$, $\tilde c
_{i}$ being a generator of $G((c _{i} \otimes _{B _{(i-1)}'} B
_{i})/B _{i})$. Then $B _{i} ^{\prime } = \widetilde \otimes _{B 
_{i}} F(r _{i}, c _{i})$, $(r _{i}, c _{i}) \in R _{p,i} \times C 
_{p,i}$; in particular, $B _{i} ^{\prime } \in {\rm Reg}(B _{(i-1)} 
^{\prime })$.
\par
\medskip
Throughout this paper, Cor$_{F/E}$ denotes the corestriction
homomorphism of Br$(F)$ into Br$(E)$, and Ker$(F/E)$ stands for the
kernel of Cor$_{F/E}$, for any finite separable extension $F/E$. The
first part of the following statement, complemented by Proposition
3.4, gives evidence of close relations between Cor$_{F/E}$ and
quasilocal nonreal fields:
\par
\medskip
(2.3) (i) $E$ is $p$-quasilocal if and only if Cor$_{R/E}$ maps
Br$(R) _{p}$ injectively into Br$(E) _{p}$ (i.e. Br$(R) _{p} \cap
{\rm Ker}(R/E) = \{0\}$), for each finite extension $R$ of $E$ in $E
(p)$ [10, (1.1) (i)];
\par
(ii) If $E$ is $p$-quasilocal and $R$ is a finite extension of $E$
in $E (p)$, then:
\par
($\alpha $) $E$ admits local $p$-CFT, provided that Br$(E) _{p} \neq
\{0\}$ [8, Theorem 3.1]; in particular, $N(R/E) \neq E ^{\ast }$
unless $R = E$;
\par
($\beta $) $N(R/E) = N(R/E) _{\rm Ab}$ [5, Theorem 3.1];
\par
($\gamma $) $N(R/E) = E ^{\ast }$ in case Br$(E) _{p} = \{0\}$ [7,
I, Lemma 4.2].
\par
\medskip
Statement (2.3) (i), Proposition 3.4 and the noted proofs of (2.3)
(ii) help us observe the possibility to apply Theorem 1.2 (i) and
the second assertion of Theorem 1.2 (ii) to the study of
Cor$_{F/E}$. This is demonstrated by the alternative proofs in [10]
of two known projection formulae (see [33, page 205]). We first
state a special case of the first projection formula, which is
particularly easy to apply.
\par
\medskip
{\bf Proposition 2.5.} {\it Let $E$, $F$ and $M$ be fields with $M
\cap F = E$, $F \subseteq E _{\rm sep}$ and $M$ cyclic over $E$, and
let $\tilde \sigma $ be a generator of $G(M/E)$. Then $MF/F$ is
cyclic, $\tilde \sigma $ extends uniquely to an $F$-automorphism
$\sigma $ of $MF$, $G(MF/F) = \langle \sigma \rangle $ and {\rm
Cor}$_{F/E}$ maps the similarity class of the cyclic $F$-algebra
$(MF/F, \sigma , \lambda )$ into $[(M/E, \tilde \sigma , N _{E} ^{F}
(\lambda ))]$, for each $\lambda \in F ^{\ast }$.}
\par
\medskip
Let now $E$ be a field and $F$ a finite extension of $E$ in $E _{\rm
sep}$, $r _{E/F}$ the restriction homomorphism $C _{E} \to C _{F}$,
and cor$_{F/E}$ the corestriction map $C _{F} \to C _{E}$. It is
known (cf. [18, Ch. 7, Corollary 5.3]) that $C _{F}$ is an abelian
torsion group and for each $p \in \overline P$, its $p$-component
can be identified with the character group $C(F (p)/F)$ of $G(F
(p)/F)$. Recall that for each $\chi \in C _{F}$, the fixed field $L
_{\chi }$ of the kernel Ker$(\chi )$ is cyclic over $F$; we denote
by $\sigma _{\chi }$ the generator of $G(L _{\chi }/F)$ induced by
any $\bar \sigma _{\chi } \in G _{E}$ satisfying the equality $\chi
(\bar \sigma _{\chi }) = (1/[L _{\chi }\colon F]) + \hbox{\Bbb Z}$.
Note that $L _{r _{E/F} (\chi )} = L _{\chi }F$, $\sigma _{r _{E/F}
(\chi )}$ is the unique $(L _{\chi } \cap F)$-automorphism of $L
_{\chi }F$ extending $\sigma ^{d(\chi )}$, and $\rho _{E/(L _{\chi }
\cap F)}$ maps Br$(L _{\chi }/E)$ on the set $\{[(L _{\chi }/(L
_{\chi } \cap F), \sigma ^{d(\chi )}, c)]\colon \ c \in E ^{\ast
}\}$, where $d(\chi ) = [L _{\chi } \cap F\colon E]$. These
observations enable one to deduce from Proposition 2.5 the first
projection formula in general (see the proof of [10, (3.1)]):
\par
\medskip
(2.4) Cor$_{F/E} ([(L _{\chi }F/F, \sigma _{r _{E/F} (\chi )},
\lambda )]) = [(L _{\chi }/E, \sigma _{\chi }, N _{E} ^{F} (\lambda
))]$, $\lambda \in F ^{\ast }$.
\par
\medskip
The second projection formula is contained in the following result,
which is used for proving Theorem 3.1:
\par
\medskip
{\bf Proposition 2.6.} {\it Let $E$ be a field, $F$ a finite
extension of $E$ in $E _{\rm sep}$, $c$ and $\chi $ elements of $E
^{\ast }$ and $C _{F}$, respectively. Then {\rm Cor}$_{F/E} ([(L
_{\chi }/F, \sigma _{\chi }, c)]) = [(L _{\tilde \chi }/E,$ $\sigma
_{\tilde \chi }, c)]$, where $\tilde \chi = {\rm cor}_{F/E} (\chi
)$. Also, $L _{\tilde \chi } \subseteq M$, provided that $M \in {\rm
Gal}(E)$ and $L _{\chi } \subseteq M$.}
\par
\medskip
The proof of Proposition 2.6 in [10] is based not only on Theorem
1.2 (i) and the second part of Theorem 1.2 (ii), applied to $T _{0}
= {\rm Br}(E _{0})$ and $\chi = {\rm Fin}$. It also relies on
Proposition 2.5, statements (2.3) (ii) ($\alpha $), ($\beta $) and
the fact (see [7, II, Lemma 2.3]) that if $F$ is $p$-quasilocal with
a primitive $p$-th root of unity, for some $p \in P(F)$, then the
structure of $C(F (p)/F)$ is determined by the group $_{p} {\rm
Br}(F) = \{b \in {\rm Br}(F)\colon \ pb = 0\}$ and the group $R _{p}
(F)$ of roots of unity in $F$ of $p$-primary degrees, as follows:
\par
\medskip
(2.5) $C(F (p)/F)$ is divisible if and only if Br$(F) _{p} = \{0\}$
or $R _{p} (F)$ is infinite. If Br$(F) _{p} \neq \{0\}$, $B _{p}$ is
a basis of $_{p} {\rm Br}(F)$ as a vector space over the field
$\hbox{\Bbb F} _{p}$ with $p$ elements, and $R _{p} (F)$ is of
finite order $p ^{\mu }$, then the group $p ^{\mu }C(F (p)/F)$ is
divisible and $C(F (p)/F)$ is isomorphic to the direct sum $p ^{\mu
}C(F (p)/F) \oplus R(F (p)/F)$, where $R(F (p)/F)$ is a subgroup of
$C(F (p)/F)$ presentable as a direct sum of cyclic groups of order
$p ^{\mu }$, indexed by $B _{p}$.
\par
\medskip
{\bf Remark 2.7.} Let $E$ be a field and $F/E$ a finite separable
extension. Suppose also that $E$ contains a primitive $n$-th root of
unity, for some $n \in \hbox{\Bbb N}$, or char$(E) = p > 0$.
Applying Kummer theory and its analogue obtained by Witt (see [21,
Ch. VIII, Sect. 8] and [18, Ch. 7, 1.9 and 2.9]), one deduces from
Proposition 2.6 the projection formula for symbol $F$-algebras of
indices dividing $n$, and the one for $p ^{\mu }$-symbol
$F$-algebras, $\mu \in \hbox{\Bbb N}$, contained in [35, Theorem
3.2] and [22, Proposition 3 (i)], respectively.
\par
\medskip
Let us mention that (2.4) and Propositions 2.5 and 2.6 can also be
proved by applying group-cohomological technique (see [40,
Proposition 4.3.7] and [18, Ch. 7, Corollary 5.3]). Without
comparing the approach referred to with the one followed in [10],
note that the latter bears an entirely field-theoretic character
both technically and conceptually. As shown in [10, Sect. 2] and
Section 5, our approach also allows us to prove Theorem 1.2 (i) and
the concluding assertion of Theorem 1.2 (ii) together with the
following result.
\par
\medskip
{\bf Proposition 2.8.} {\it Let $F$ be a field, $M/F$ a finite
separable extension, $\Lambda ^{\prime }$ a tensor compositum over
$M$ of function fields of Brauer-Severi $M$-varieties, and $\Lambda
$ is the transfer of $\Lambda ^{\prime }$ over $M/F$. Then {\rm
Br}$(\Lambda /F)$ equals the image of {\rm Br}$(\Lambda ^{\prime
}/M)$ under {\rm Cor}$_{M/F}$.}
\par
\medskip
Proposition 2.8 is a special case of [13, Proposition 2.6] which has
been deduced in [13] from [30, Theorem 3.13] and the description of
the relative Brauer groups of function fields of generalized
Brauer-Severi varieties [3]. In view of the relations between
quasilocal nonreal fields and Brauer group corestrictions, and of
the preservation of rationality under transfer (see [30, Lemma 3.2
(a)]), one may expect that Theorem 1.2 can be used for simplifying
the proofs and the presentations of index reduction formulae, for
the function fields of a number of twisted rational varieties like
those considered in [30; 24; 23] and [37].
\par
\vskip0.6truecm \centerline {\bf 3. $p$-primary analogue to Theorem
1.1}
\par
\medskip
Let $E$ be a field, $R$ a finite extension of $E$ in $E _{\rm sep}$,
and $H(E) ^{n} = \{h ^{n}\colon \ h \in H(E)\}$, for any subgroup
$H(E)$ of $E ^{\ast }$ and each $n \in \hbox{\Bbb N}$. For each $p
\in \overline P$, let $R _{{\rm Ab},p} = R \cap E _{\rm Ab} \cap E
(p)$, $\rho _{p}$ be the greatest divisor of $[R\colon E]$, and $N
_{p} (R/E) = \{u _{p} \in E ^{\ast }\colon $ the co-set $u
_{p}N(R/E)$ is a $p$-element of $E ^{\ast }/N(R/E)\}$. Clearly, $E
^{\ast \rho _{p}} \subseteq N _{p} (R/E)$ and $N(R/E) _{\rm Ab}
^{\rho _{p}} \subseteq N(R _{{\rm Ab},p}/E)$, $p \in \overline P$,
so Theorem 1.1 can be deduced from the following result (in the case
of $\Omega = E _{\rm sep}$):
\par
\medskip
{\bf Theorem 3.1.} {\it Let $E$ be a field, $\Omega $ a Galois
extension of $E$ in $E _{\rm sep}$, and $R$ a finite subextension of
$E$ in $\Omega $. Assume that finite extensions of $E$ in $\Omega $
are $p$-quasilocal, for some $p \in \overline P$. Then $N(R/E) =
N(R _{{\rm Ab},p}/E) \cap N _{p} (R/E)$ in the following cases:}
\par
(i) {\it There exists $M \in {\rm Gal}(E) \cap I(\Omega /E)$ with $R
\in I(M/E)$ and {\rm Br}$(M) _{p} \subseteq {\rm Im}(E/M)$;}
\par
(ii) $N(R/E)$ {\it includes $N(\Phi (R)/E) \cap N _{p} (R/E)$, for
some $\Phi (R) \in \Omega _{p} (E) \cap I(\Omega /E)$.}
\par
\medskip
{\it Proof.} The inclusion $N(R/E) \subseteq N(R _{{\rm Ab},p}/E)
\cap N _{p} (R/E)$ is obvious. We prove the converse by showing
that $p \not\vert e _{p}$, where $e _{p}$ is the exponent of the
$N(R _{{\rm Ab},p}/E)/N(R/E)$. Theorem 3.1 (ii) is obtained as a
special case of the following lemma.
\par
\medskip
{\bf Lemma 3.2.} {\it Let $E$, $\Omega $, $R$ and $p$ satisfy the
conditions of Theorem 3.1, and let $M \in {\rm Gal}(E) \cap
I(\Omega /E)$. Then $N _{p} (M/E)N _{p} (R/E) = N _{p} (M \cap
R/E)$.}
\par
\medskip
{\it Proof.} It is clearly sufficient to show that $p$ does not
divide the exponent $e$ of $L ^{\ast }/N(M/L)N(R/L)$, where $L = M
\cap R$. Hence, by the $p$-quasilocal property of finite extensions
of $L$, one may assume further that $L = E$. Let $\widetilde R$ be
the normal closure of $R$ in $E _{\rm sep}$ over $E$, $\widetilde H
_{p} \in {\rm Syl}_{p} (\widetilde R/R)$, $\widetilde G _{p} \in
{\rm Syl}_{p} (\widetilde R/E)$, $\widetilde H _{p} \subseteq
\widetilde G _{p}$ and $G _{p} \in {\rm Syl}_{p} (M/E)$. Denote by
$R _{p}$, $\Phi _{p}$ and $E _{p}$ the fixed fields of $\widetilde H
_{p}$, $\widetilde G _{p}$ and $G _{p}$, respectively. We prove that
$e \vert \varphi _{p}$, i.e. $E ^{\ast \varphi _{p}} \subseteq
N(M/E)N(R/E)$, where $\varphi _{p} = [\Phi _{p}\colon E][E
_{p}\colon E]$; this implies the lemma, since $p \not\vert \varphi
_{p}$. Let $\xi $ be a primitive element of $R/E$. By Galois theory
and the equality $L = E$, $[MR\colon R] = [M\colon E]$, i.e.
$[M(\xi )\colon M] = [R\colon E]$. This means that the minimal
polynomial of $\xi $ over $E$ is irreducible over $M$. Considering
$E _{p}$, $RE _{p}$ and $M$ instead of $E$, $R$ and $M$,
respectively, and using Lemma 2.2, one reduces our proof to the
special case where $E _{p} = E$, i.e. $M \subseteq E (p)$. The
choice of $R _{p}$ guarantees that $p \not\vert [R _{p}\colon R]$,
so it follows from Galois theory, the equality $[M\colon E] =
[MR\colon R]$ and the inclusion $M \subseteq E (p)$ that $[MR
_{p}\colon R _{p}] = [M\colon E]$ and $R _{p} \cap MF = F$, for
every $F \in I(R _{p}/E)$. This, applied to the case of $F = \Phi
_{p}$, enables one to deduce from (2.3) (ii) ($\alpha $) and
($\gamma $) that $N(M\Phi _{p}/\Phi _{p})N(R _{p}/\Phi _{p}) =
\Phi _{p} ^{\ast }$. Now the inclusion $E ^{\ast \varphi _{p}}
\subseteq N(M/E)N (R/E)$ becomes obvious, so Lemma 3.2 is proved.
\par
\medskip
{\bf Remark 3.3.} Statements (1.2) (i) and (ii) show that the
conclusions of Theorem 1.1 (ii) and Lemma 3.2 are not always true,
if $E$ is only a PQL-field.
\par
\medskip
Our next result characterizes the fields whose finite extensions are
$p$-quasilocal, for a given $p \in \overline P$. It simplifies the
proof of Theorem 3.1 (i) and leads to the idea of constructing
quasilocal nonreal fields and formally real PQL-fields by the method
followed in this paper.
\par
\medskip
{\bf Proposition 3.4.} {\it Let $E$ be a field and $\Omega $ a
Galois extension of $E$ in $E _{\rm sep}$. Then finite extensions of
$E$ in $\Omega $ are $p$-quasilocal, for a given $p \in \overline
P$, if and only if one of the following two conditions is fulfilled:}
\par
(c) $p > 2$ {\it or $E$ is nonreal, and for each pair $(M, M
^{\prime }) \in {\rm Gal}(E) \times {\rm Gal}(E)$ with $M \subseteq
\Omega $ and $M ^{\prime } \in I(M (p)/M)$, ${\rm Cor}_{M'/M}$ maps
{\rm Br}$(M ^{\prime }) _{p}$ injectively into {\rm Br}$(M) _{p}$;}
\par
(cc) $p = 2${\it , $E$ is formally real and its formally real
finite extensions in $\Omega $ are uniquely
ordered Pythagorean fields; when $\Omega \neq E$, this holds if and
only if $2 \not\in P(\Omega )$, $G(\Omega /E(\sqrt{-1}))$ is
abelian, {\rm cd}$_{2} (G(\Omega /E(\sqrt{-1})) = 0$ and $G(\Omega
/E)$ is continuously isomorphic to the semidirect product $G(\Omega
/E(\sqrt{-1})) \times \langle \sigma \rangle $, where $\sigma ^{2} =
1$ and $\sigma \tau \sigma = \tau ^{-1}$, for all $\tau
\in G(\Omega /E(\sqrt{-1}))$.
\par
If (cc) is satisfied and $\Omega = E _{\rm sep}$, then $N(L/E) = N
_{2} (L/E) \cap N(L/E) _{\rm Ab}$, {\rm Br}$(E)$ is of order $2$ and
{\rm Im}$(E/L) = {\rm Br}(L) _{2}$, for each finite extension $L$ of
$E$ in $E _{\rm sep}$.}
\par
\medskip
{\it Proof.} It is clear from (2.3) (i) that if finite extensions of
$E$ in $\Omega $ are $p$-quasilocal, then Cor$_{F'/F}$ maps Br$(F
^{\prime }) _{p}$ injectively into Br$(F) _{p}$ whenever $F \in {\rm
Gal}(E)$, $F \subseteq \Omega $ and $F ^{\prime }$ is a finite
extension of $F$ in $F (p)$. We first prove that the fulfillment of
(c) implies that finite extensions of $E$ in $\Omega $ are
$p$-quasilocal. To begin with, (2.3) (i) and Proposition 2.3 (ii)
guarantee that $E$ is $p$-quasilocal and Br$(E) _{p}$ is divisible.
Observe that both properties are preserved by each $M \in {\rm
Gal}(E)$, $M \subseteq \Omega $. If $p \not\in P(M)$, this follows
at once from Proposition 2.3 (ii), so we assume further that $p \in
P(M)$. Denote by $\hbox{\Bbb F}$ the prime subfield of $E$, and by
$\Gamma _{p}$ the unique $\hbox{\Bbb Z} _{p}$-extension of
$\hbox{\Bbb F}$ in $E _{\rm sep}$. The divisibility of Br$(M) _{p}$
can be deduced from Witt's theorem (see [11, Sect. 15]), if $p =
{\rm char}(E)$, and from the Merkurjev-Suslin theorem [25, (16.1)]
in case $p \neq {\rm char}(E)$ and $\Gamma _{p} \subseteq E$.
Assuming that $\Gamma _{p} \not\subseteq E$, one obtains from Galois
theory that $M\Gamma _{p}/M$ is a $\hbox{\Bbb Z} _{p}$-extension, and
for each $n \in \hbox{\Bbb N}$, Gal$(E)$ contains the field $M _{n}
\in I(M\Gamma _{p}/M)$ of degree $[M _{n}\colon M] = p ^{n}$. Hence,
by (c) and the RC-formula, $M _{n}$ is a splitting field of each $T
_{n} \in d(M _{n})$ of exponent dividing $p ^{n}$. This enables one
to deduce the following statements, arguing as in the proof of [7,
I, Theorem 3.1]:
\par
\medskip
(3.1) (i) ind$(\Delta ) = {\rm exp}(\Delta )$, for every $\Delta \in
d(M)$ of $p$-primary dimension.
\par
(ii) If Br$(M) _{p} \neq \{0\}$, then there exists $\Delta _{n} \in
d(M)$ of index $p ^{n}$, for each $n \in \hbox{\Bbb N}$.
\par
\medskip
Statement (3.1) the established property of $M _{n}$ and [27, Sect.
15.1, Corollary b] imply the divisibility of Br$(M) _{p}$. Our
objective now is to prove that $M$ is $p$-quasilocal, provided that
$p \in P(M)$. Let $\widetilde M$ be a finite extension of $M$ in $M
(p)$, $D$ a central division $M$-algebra, such that exp$(D) = [M
^{\prime }\colon M] = p$. Then it follows from Galois theory that $M
(p)$ contains as a subfield the normal closure $\widetilde M
^{\prime }$ of $\widetilde M$ over $E$ (in $E _{\rm sep}$). Since
Br$(M) _{p}$ is divisible, one can find an algebra $\widetilde D \in
d(M)$ so as to satisfy the equalities ind$(\widetilde D) = p
^{\tilde m}$ and $p ^{\tilde m-1}\widetilde [D] = [D]$, where $p
^{\tilde m} = [\widetilde M ^{\prime }\colon M]$. Hence, by the
RC-formula, Cor$_{\widetilde M'/\widetilde M} (\rho _{M/\widetilde
M'} ([\widetilde D]) = [D]$, and Cor$_{\widetilde M'/M} ([\widetilde
D]) = 0 = {\rm Cor}_{\widetilde M/M} (\rho _{M/\widetilde M} ([D])$.
In view of (c), this means that $\rho _{M/\widetilde M'}
([\widetilde D]) = 0$ and $\rho _{M/\widetilde M} ([D]) = 0$ (in
Br$(\widetilde M ^{\prime })$ and Br$(\widetilde M)$, respectively).
In other words, $\widetilde D$ is split by $\widetilde M ^{\prime }$
and $D$ is split by $\widetilde M$, so the $p$-quasilocal property
of $M$ becomes obvious.
\par
Suppose now that $R$ is an arbitrary finite extension of $E$ in
$\Omega $ and denote by $R _{1}$ its normal closure in $E _{\rm
sep}$ over $E$. By definition, $R$ is $p$-quasilocal, if $p \not\in
P(R)$ or Br$(R) _{p} = \{0\}$, so we assume that $p \in P(R)$ and
Br$(R) _{p} \neq \{0\}$. Note first that Br$(R) _{p}$ is divisible.
As in the special case where $R = R _{1}$, one sees that it is
sufficient to prove our assertion under the hypothesis that $\Gamma
_{p} \not\subseteq E$. Applying [27, Sect. 15.1, Corollary b], one
concludes that if Br$(R) _{p}$ is not divisible, then Br$(R\Gamma
_{p}/R) \neq {\rm Br}(R) _{p}$ and Br$(R\Gamma _{p}) _{p}$ is
infinite and divisible. As $[R _{1}\Gamma _{p}\colon R\Gamma _{p}]
\in \hbox{\Bbb N}$, this implies that Br$(R _{1}\Gamma _{p}) _{p}
\neq \{0\}$. On the other hand, $R _{1}$ is $p$-quasilocal and $R
_{1}\Gamma _{p}/R _{1}$ is a $\hbox{\Bbb Z} _{p}$-extension, so it
follows from [7, I, Theorem 4.1 (iv)] that Br$(R _{1}\Gamma _{p})
_{p} = \{0\}$. The obtained contradiction completes the proof of the
divisibility of Br$(R) _{p}$, so we return to the assumption that $p
\in P(R)$ and Br$(R) _{p} \neq \{0\}$. Let $R ^{\prime }$ be an
extension of $R$ in $R (p)$ of degree $p$, $R _{1} ^{\prime }$ the
normal closure of $R ^{\prime }$ in $E _{\rm sep}$ over $E$, and $d
_{1}$ an element of Br$(R ^{\prime }) _{p}$, such that Cor$_{R'/R}
(d _{1}) = 0$. Then Br$(R ^{\prime } _{p})$ is divisible, so the
equation $[R _{1} ^{\prime }\colon R ^{\prime }]x = d _{1}$ has a
solution $d _{1} ^{\prime } \in {\rm Br}(R ^{\prime }) _{p}$. Hence,
by the RC-formula, Cor$_{R _{1}'/R'} (d _{1} ^{\prime }) = d _{1}$.
In view of the equality Cor$_{R _{1}'/R} = {\rm Cor}_{R'/R} \circ
{\rm Cor}_{R _{1}'/R'}$, this means that $\rho _{R'/R _{1}'} (d _{1}
^{\prime }) \in {\rm Ker}(R _{1}'/R)$. Thus condition (c) yields
$\rho _{R/R _{1}'} (d _{1}) = 0$ and $d _{1} = 0$. The obtained
result indicates that Cor$_{R'/R}$ maps Br$(R ^{\prime }) _{p}$
injectively into Br$(R) _{p}$, so the assertion that $R$ is
$p$-quasilocal reduces to a consequence of (2.3) (i).
\par
Assume now that $E$ is formally real and $p = 2$. Note first that it
suffices for the proof of Proposition 3.4 (cc) to show that if
$\Omega \neq E$ and formally real finite extensions of $E$ in
$\Omega $ are $2$-quasilocal, then $2 \not\in P(\Omega )$ and
Br$(\Phi ) _{2} = \{0\}$, for each $\Phi \in {\rm Gal}(E) \cap
I(\Omega /E)$, $\Phi \neq E$. This follows from Becker's theorem
(cf. [4, (3.3)]), [9, Proposition 3.1] and the latter part of
Proposition 2.3 (ii). Observe that every admissible $\Phi $ is a
nonreal field. Indeed, for each primitive element $\xi $ of $\Phi
/E$ and any $g \in G(\Phi /E)$, the trace Tr$_{E} ^{\Phi } (\xi -
g(\xi ))$ equals zero. Therefore, the hypothesis that $\Phi $ is
formally real requires that $(\xi - g(\xi )) \not\in \Phi ^{\ast 2}
\cup (-1)\Phi ^{\ast 2}$. In view of Proposition 2.3 (ii), when $g
\neq 1$, this contradicts the assumption that $\Phi $ is
$2$-quasilocal, so the assertion that $\Phi $ is nonreal is proved.
Let $E _{2}$ be the fixed field of some $G _{2} \in {\rm Syl}_{2}
(\Phi /E)$. Then $2 \not\vert [E _{2}\colon E]$, and by the
Artin-Schreier theory (cf. [21, Ch. XI, Proposition 2]), $E _{2}$
is formally real. Hence, by Proposition 2.3 (ii), $\Phi = E _{2}
(\sqrt{-1})$, $2 \not\in P(\Phi )$ and $2 \not\vert [\Phi \colon
E(\sqrt{-1})]$, so Br$(\Phi ) _{2} = \{0\}$ and $2 \not\in P(\Omega
)$, which proves (cc).
\par
Henceforth, we assume that $\Omega = E _{\rm sep}$. This ensures
that $P(E) = \{2\}$, $E _{\rm Ab} = E(\sqrt{-1})$ and Br$(E)$ is of
order $2$ (cf. [4, (3.3)] and [9, Lemma 2.4]). Let $L$ be a finite
extension of $E$ in $E _{\rm sep}$. By [9, Proposition 3.1], $L \in
{\rm Gal}(E)$ if and only if $\sqrt{-1} \in L$. When $L \in {\rm
Gal}(E)$, the same result shows that $2 \not\vert [L\colon
E(\sqrt{-1})]$, which yields Br$(L) _{2} = \{0\}$ and $N _{2}
(L/E(\sqrt{-1})) = E(\sqrt{-1}) ^{\ast }$. Since $E ^{\ast } = E
^{\ast 2} \cup (-1)E ^{\ast 2}$ and, by Lemma 2.2, $N _{2} (L/E)
\subseteq N _{2} (L/E(\sqrt{-1}))$, this means that $N _{2} (L/E)
\cap N(L/E) _{\rm Ab} = N _{2} (L/E) ^{2}$ $= N(L/E)$. Suppose
finally that $L \not\in {\rm Gal}(E)$. Then $2 \not\vert [L\colon
E]$ and $L$ is formally real, which implies that $\rho _{E/L}$ is an
isomorphism and $N _{2} (L/E) = N(L/E)$. These results prove
Proposition 3.4.
\par
\medskip
Theorem 3.1 (ii), Propositions 2.3 (ii)-3.4 and our next statement
reduce the proof of Theorem 3.1 (i) to the case where Br$(E) _{p}$
is divisible and $R = M$.
\par
\medskip
{\bf Corollary 3.5.} {\it Assume that $E$ and $M$ are fields, such
that $M \in {\rm Gal}(E)$, {\rm Br}$(E) _{p}$ is divisible and
$I(M/E)$ consists of $p$-quasilocal fields, for some $p \in
\overline P$. Then {\rm Br}$(L) _{p}$ is divisible, for every
extension $L$ of $E$ in $M$, and the following conditions are
equivalent:}
\par
(i) $\rho _{E/M}$ {\it maps {\rm Br}$(E) _{p}$ surjectively on {\rm
Br}$(M) _{p}$;}
\par
(ii) Cor$_{M/E}$ {\it maps {\rm Br}$(M) _{p}$ injectively into {\rm
Br}$(E) _{p}$;}
\par
(iii) {\it For each $R \in I(M/E)$, {\rm Br}$(R) _{p} \subseteq {\rm
Im}(E/R)$ and {\rm Br}$(R) _{p} \cap {\rm Ker}(R/E) = \{0\}$.}
\par
\medskip
{\it Proof.} The conclusion that Br$(L) _{p}$ is divisible follows
from Proposition 2.3 (i), the divisibility of Br$(E) _{p}$ and the
$p$-quasilocal property of $E$ and $L$. For any pair $(U, V) \in
I(M/E) \times I(M/E)$, such that $U \subseteq V$, and put Im$
(U/V) _{p} = {\rm Br}(V) _{p} \cap {\rm Im}(U/V)$ and Ker$(V/U) _{p}
= {\rm Br}(V) _{p} \cap {\rm Ker}(V/U)$. Since Br$(U) _{p}$ is
divisible, the RC-formula implies Br$(V) _{p} = {\rm Im}(U/V) _{p} +
{\rm Ker}(V/U) _{p}$. This, applied to $(E, M)$, proves that
(ii)$\to $(i).
\par
The rest of the proof relies on the well-known fact (see, e.g. [35])
that for each tower $U _{1} \subseteq U _{2} \subseteq U _{3}$ of
finite separable extensions, Cor$_{U _{3}/U _{1}}$ equals the
composition ${\rm Cor}_{U _{2}/U _{1}} \circ {\rm Cor}_{U _{3}/U
_{2}}$. This implies that if $M/E$ satisfies (ii), then so does
$M/U$. At the same time, it follows from the RC-formula and the
divisibility of Br$(V) _{p}$ that Br$(V) _{p}$ is included in the
image of Br$(V) _{p}$ under Cor$_{M/V}$. Considering now the tower
$U \subseteq V \subseteq M$, one concludes that if condition (ii)
holds, then Cor$_{V/U}$ maps Br$(V) _{p}$ injectively into Br$(U)
_{p}$.
\par
Let now $E _{p}$ be the fixed field of some $G _{p} \in {\rm
Syl}_{p}(M/E)$. Then $p \not\vert [E _{p}\colon E]$, so the
RC-formula and the general properties of Schur indices (cf. [27,
Sect. 13.4]) imply that the sum Br$(E _{p}) _{p} = {\rm Im}(E/E
_{p}) _{p} + {\rm Ker}(E _{p}/E) _{p}$ is direct. Since $\rho _{E/M}
= \rho _{E _{p}/M} \circ \rho _{E/E _{p}}$, one also sees that if
condition (i) holds, then Br$(E _{p}) _{p} = {\rm Im}(E/E _{p}) _{p}
+ {\rm Br}(E _{p}/M)$. As Br$(E _{p}) _{p}$ is divisible and Br$(E
_{p}/M)$ is of exponent dividing $[M\colon E _{p}]$, these
observations yield Br$(E _{p}) _{p} = {\rm Im}(E/E _{p}) _{p}$ and
Ker$(E _{p}/E) _{p} = \{0\}$. It is now easily obtained from (2.3)
(i), applied to $E _{p}$ and $M$, and the equality Cor$_{M/E} = {\rm
Cor}_{E _{p}/E} \circ {\rm Cor}_{M/E _{p}}$ that (i)$\to $(ii).
Returning to the beginning of our proof, one also sees that (i)
implies Br$(M) _{p} = {\rm Im}(U/M) _{p}$ and Br$(V) _{p} \subseteq
{\rm Im}(U/V) _{p} + {\rm Br}(M/V)$. Since Br$(V) _{p}$ is divisible
and the exponent of Br$(M/V)$ divides $[M\colon V]$, this yields
Br$(V) _{p} = {\rm Im}(U/V) _{p}$, which completes the proof of the
implication (i)$\to $(iii). As (iii) obviously implies (i),
Corollary 3.5 is proved.
\par
\medskip
Now we prove Theorem 3.1 (i) in the case where Br$(E) _{p} = \{0\}$
and $R = M$. Let $E _{p}$ be the fixed field of a group $G _{p} \in
{\rm Syl}_{p} (M/E)$, and let $[E _{p}\colon E] = m _{p}$. By
Corollary 3.5, Br$(E _{p}) _{p} = \{0\}$, so it follows from (2.3)
(ii) ($\gamma $) that $N(M/E _{p}) = E _{p} ^{\ast }$. Hence, by the
norm equality $N _{E} ^{M} = N _{E} ^{E _{p}} \circ N _{E _{p}}
^{M}$, $N(M/E) = N(E _{p}/E)$. As $p \not\vert m _{p}$ and $E ^{\ast
m _{p}} \subseteq N(E _{p}/E)$, this implies that $N _{p} (M/E)
\subseteq N(M/E)$, as claimed. For the proof of Theorem 3.1 (i) in
the case of Br$(E) _{p} \neq \{0\}$, we need the following lemmas.
\par
\medskip
{\bf Lemma 3.6.} {\it Let $E$ be a field, $M \in {\rm Gal}(E)$, $G
_{p} \in {\rm Syl}_{p} (M/E)$, for some $p \in \overline P$, $E
_{p}$ the fixed field of $G _{p}$, $L$ a cyclic extension of $E
_{p}$ in $M$, and $\sigma $ a generator of $G(M/E)$. Assume that $p$
does not divide the index of the commutator subgroup $[G(M/E),
G(M/E)]$ in $G(M/E)$. Then $[(L/E _{p}, \sigma , c)] \in {\rm Ker}(E
_{p}/E)$, for every $c \in E ^{\ast }$.}
\par
\medskip
{\it Proof.} Our assumptions show that $M \cap E (p) = E$. Since $L
= L _{\chi }$, for some $\chi \in C(F (p)/F)$, this reduces our
assertion to a consequence of Proposition 2.6.
\par
\medskip
{\bf Lemma 3.7.} {\it With assumptions being as in Theorem 3.1, let
$M \cap E (p) = E$. Then $E ^{\ast } \subseteq N(M/E _{p})$ and $N
_{p} (R/E) \subseteq N(R/E)$.}
\par
\medskip
{\it Proof.} Clearly, one may consider only the special case of $R =
M \neq E$ and Br$(E) _{p} \neq \{0\}$. Take $G _{p}$ and $E _{p}$ as
in Lemma 3.6 and put $m _{p} = [E _{p}\colon E]$. We show that $E
^{\ast } \subseteq N(M/E _{p})$. As $p \not\vert m _{p}$, $\rho
_{E/E _{p}}$ maps Br$(E) _{p}$ injectively into Br$(E _{p}) _{p}$.
Therefore, Br$(E _{p}) _{p} \neq \{0\}$ and since $E _{p}$ is
$p$-quasilocal, (2.3) (ii) ($\alpha $) and ($\beta $) indicate that
it is sufficient to prove the inclusion $E ^{\ast } \subseteq N(L/E
_{p})$, for an arbitrary cyclic extension $L$ of $E _{p}$ in $M$. By
[27, Sect. 15.1, Proposition b], this amounts to showing that $[(L/E
_{p}, \sigma , c)] = 0$ in Br$(E _{p})$, for each $c \in E ^{\ast
}$, where $\sigma $ is a generator of $G(L/E _{p})$. As g.c.d.$(m
_{p}, p) = 1$ and $\rho _{E/E _{p}}$ maps Br$(E) _{p}$ surjectively
on Br$(E _{p}) _{p}$, Cor$_{E _{p}/E}$ induces an isomorphism Br$(E
_{p}) _{p} \cong {\rm Br}(E) _{p}$ (cf. [35, Theorem 2.5]). This,
combined with Lemma 3.6, implies that $[(L/E _{p}, \sigma , c)] =
0$, $c \in E ^{\ast }$, as claimed. Hence, $E ^{\ast } \subseteq
N(M/E _{p})$ and $E ^{\ast m _{p}} \subseteq N(M/E)$, which proves
Lemma 3.7.
\par
\medskip
{\bf Lemma 3.8.} {\it In the setting of Theorem 3.1 (i), let
$M \in {\rm Gal}(E)$ and $G(M/E) \in {\rm Sol}$. Then $N _{p} (M/E)
\cap N(M _{{\rm Ab},p}/E) \subseteq N(M/E)$.}
\par
\medskip
{\it Proof.} It is sufficient to prove the lemma under the
hypothesis that $N(M ^{\prime }/E ^{\prime })$ includes $N _{p} (M
^{\prime }/E ^{\prime }) \cap N(M ^{\prime } _{{\rm Ab},p}/E
^{\prime })$, whenever $E ^{\prime }$ and $p$ satisfy the conditions
of Theorem 3.1, $M ^{\prime } \in {\rm Gal}(E ^{\prime })$, $G(M
^{\prime }/E ^{\prime }) \in {\rm Sol}$ and $[M ^{\prime }\colon E
^{\prime }] < [M\colon E]$. As in the proof of [5, Theorem 1.1], we
first show that one may assume further that $G(M/E)$ is a
Miller-Moreno group (i.e. nonabelian whose proper subgroups lie in
${\rm Ab}$). Our argument relies on the fact that the class of
fields satisfying the conditions of Theorem 3.1 is closed under the
formation of finite extensions. Note that if $G(M/E)$ is not
Miller-Moreno, then it has a subgroup $H \not\in {\rm Ab}$ with $[H,
H]$ normal in $G(M/E)$. Indeed, one can put $H = [G(M/E), G(M/E)]$
in case $G(M/E) \not\in {\rm Met}$, and take as $H$ any nonabelian
maximal subgroup of $G(M/E)$, otherwise. Let $F$ and $L$ be the
fixed fields of $H$ and $[H, H]$, respectively. It follows from
Galois theory and the choice of $H$ that $L \in {\rm Gal}(E)$, $M
\cap E _{\rm Ab} \subseteq L$ and $E \neq L \neq M$, so our extra
assumption and Lemma 2.1 lead to the conclusion that $N _{p} (L/E)
\cap N(M _{{\rm Ab},p}/E) \subseteq N(L/E)$ and $N _{p} (M/F) \cap
N(L/F) \subseteq N(M/F)$. Let now $\mu $ be an element of $N _{p}
(M/E) \cap N(M _{{\rm Ab},p}/E)$, and $\lambda \in L ^{\ast }$ be of
norm $N _{E} ^{L} (\lambda ) = \mu $. Then $N _{F} ^{L} (\lambda )
^{k} \in N _{p} (M/F)$, for some $k \in \hbox{\Bbb Z}$ such that $p
\not\vert \ k$. Therefore, $N _{F} ^{L} (\lambda ) ^{k} \in N(M/F)$
and $\mu ^{k} \in N(M/E)$. As $\mu \in N _{p} (M/E)$, this implies
that $\mu \in N(M/E)$, which gives the desired reduction. In view of
(2.3) (ii) ($\beta $), Lemma 2.1 and Galois theory, one may assume
that $G(M/E)$ is a Miller-Moreno group and $G(M/E) \not\in {\rm
Nil}$. Denote by $\theta $ the order of $[G(M/E), G(M/E)]$. The
assertion of the lemma is obvious, if $p \not\vert \ \theta $, so we
suppose further that $p \vert \theta $. By Miller-Moreno's
classification of these groups or by Schmidt's theorem (cf. [29,
Theorem 445] and [34, Theorem 26.1]), $G(M/E)$ has the following
properties:
\par
\medskip
(3.2) $[G(M/E), G(M/E)]$ is a minimal normal subgroup of
$G(M/E)$, which lies in Syl$_{p} (M/E) \cap {\rm Ab}$ and has
exponent $p$. Also, $G(M/E)/([G(M/E),$ $G(M/E)]) \cong C _{\pi
^{n}}$, with $C _{\pi ^{n}} \subset G(M/E)$ cyclic of order $\pi
^{n}$, for some $\pi \in \overline P$ and $n \in \hbox{\Bbb N}$. The
centre of $G(M/E)$ equals the subgroup of $C _{\pi ^{n}}$ of order
$\pi ^{n-1}$, and $\theta = p ^{k}$, where $k$ is the order of $p$
modulo $\pi $.
\par
\medskip
It follows from (3.2) and Galois theory that $[E _{p}\colon E] = \pi
^{n}$, where $E _{p} = M \cap E _{\rm Ab}$. Hence, by Lemma 3.7, $E
^{\ast } \subseteq N(M/E _{p})$ and $E ^{\ast \pi ^{n}} \subseteq
N(M/E)$. Since $N _{E _{p}} ^{M} (\eta ) = \eta ^{p ^{k}}$, for
every $\eta \in E _{p}$, we also have $N(M/E) _{\rm Ab} ^{p ^{k}}
\subseteq N(M/E)$, so Lemma 3.8 is proved.
\par
\medskip
It is now easy to complete the proof of Theorem 3.1 (i). Put $M _{0}
= M \cap E _{\rm Sol}$, and denote by $\mu _{p}$ and $m _{p}$ the
maximal divisors of $[M _{0}\colon E]$ and $[M\colon E]$,
respectively, for which $p \not\vert \mu _{p}m _{p}$. Using
Corollary 3.5 as well as the inclusion $(M \cap E _{\rm Ab})
\subseteq M _{0}$, and applying Lemma 3.8 to $M _{0}/E$ and Lemma
3.7 to $M/M _{0}$, one obtains that $N(M/E) _{\rm Ab} ^{\mu _{p}}
\subseteq N(M _{0}/E)$ and $M _{0} ^{\ast \bar m _{p}} \subseteq
N(M/M _{0})$, where $\bar m _{p} = m _{p}/\mu _{p}$. Hence, by the
norm identity $N _{E} ^{M} = N _{E} ^{M _{0}} \circ N _{M _{0}}
^{M}$, we have $N(M/E) _{\rm Ab} ^{m _{p}} \subseteq N(M _{0}/E)
^{\bar m _{p}} \subseteq N(M/E)$, i.e. $m _{p}$ is divisible by the
exponent of $N(M/E) _{\rm Ab}/N(M/E)$, so Theorem 3.1 (i) is proved.
\par
\medskip
{\bf Remark 3.9.} Theorem 3.1 (i) plays a role in the proof of the
first of the following two results (see the references in [6, page
384]):
\par
(i) There exists a nonreal SQL-field $E$, such that $G _{E}$
is not pronilpotent but is metabelian and every finite extension
$R$ of $E$ in $E _{\rm sep}$ is subject to the alternative $R
\subseteq E _{\rm Nil}$ or $N(R/E) \neq N(\Theta /E)$, $\Theta \in
\Omega (E)$.
\par
(ii) Let $F$ be a formally real quasilocal field and $\Phi $ a finite
extension of $F$ in $F _{\rm sep}$. We have already proved that then
$G _{F}$ is metabelian, $F _{\rm Nil} = F(\sqrt{-1})$,
$\Phi (\sqrt{-1}) \in {\rm Gal}(F)$, and $\rho _{F/\Phi _{0}}$ is
surjective, for each formally real field $\Phi _{0} \in I(\Phi /F)$.
Note also that the following conditions are equivalent: (a) $N(\Phi
/F) = N(\Phi /F) _{\rm Ab}$; (b) $\rho _{F/\Phi (\sqrt{-1})}$ is
surjective; (c) cd$_{t} (G _{F}) = 1\colon $ $t \in \overline P
\setminus \{2\}$, $t \vert [\Phi \colon F]$; when $F$ is not SQL, 
this holds in infinitely many cases (see [9, Lemma 2.3 and Remark 
3.2]). On the contrary, it follows from [9, Theorems 1.1 and 1.2]
that if $F$ is SQL, then $N(\Phi /F)$ is uniquely determined by the
$F$-isomorphism class of $\Phi $; hence, $\Phi $ is subject to the
alternative in (i) unless $F$ is real closed.
\par
\medskip
{\bf Corollary 3.10.} {\it Let $E$ be a quasilocal field, and
suppose that $M \in {\rm Gal}(E)$ has the property that $\rho
_{E/L}$ is surjective, where $L$ is the fixed field of the Fitting
subgroup of $G(M/E)$. Then $N(R/E) = N(R/E) _{\rm Ab}$, for each $R
\in I(M/E)$.}
\par
\medskip
{\it Proof.} In view of Theorem 1.1 (i), one may consider the
special case where $M \neq L$. It is easily seen that the exponent
of $N(L/E)/N(M/E)$ divides $[M\colon L]$. Observe also that, by
Fitting's theorem, $G(M/L)$ is normal in $G(M/E)$ and $G(M/L) \in
{\rm Nil}$. Hence, by Galois theory, the field $M \cap L (p) \colon
= M _{p}$ lies in Gal$(E) \cap I(L (p)/L)$, and by Lemma 2.1, $p
\not\vert [M\colon M _{p}]$, for any $p \in \overline P$, $p \vert
[M\colon L]$. This shows that $p$ does not divide the exponents of
$N(M/M _{p})$ and $N(M/E)/N(M _{p}/E)$. At the same time, by
Proposition 2.3 (iii), $M _{p}/E$ and $p$ satisfy the conditions of
Theorem 3.1 (i). Since $p \not\vert [M\colon M _{p}]$, one deduces
from Galois theory that $M _{{\rm Ab},p} \subseteq M _{p}$, so the
obtained results imply $p$ does not divide the exponent of
$N(M/E)/N(M/E) _{\rm Ab}$ either. Thus it follows that $N(R/E) =
N(R/E) _{\rm Ab}$, as claimed.
\par
\medskip
It follows from Corollary 3.5 and Remark 3.9 (ii) that the
conditions of Theorem 1.1 (i) guarantee the surjectivity of $\rho
_{E/R}$, for all $R \in I(M/E)$. This means that Theorem 1.1 (i) is
a special case of Corollary 3.10. Remark 3.9 (ii) shows that the
conclusion of Theorem 1.1 (i) is not necessarily true without the
assumption that $M \in {\rm Gal}(E)$. We prove in Section 6 that
the scope of Corollary 3.10 is larger than that of Theorem 1.1 (i).
\par
\vskip0.6truecm \centerline{\bf 4. On the relative Brauer group of
function field extensions}
\smallskip
\centerline{\bf of arbitrary fields, associated with norm
equations}
\par
\medskip
The results of this Section form the technical basis for the proof
of Theorem 1.2 (iii) and (iv). The main one is obtained by the
method of proving Proposition 2.6, using at crucial points (2.3)
(ii), Theorem 1.1 (i), the regularity and other known properties of
function fields of Brauer-Severi varieties and of their transfers
over finite Galois extensions (see [14, Theorem 1] with its proof).
It illustrates the fact that the applications of $p$-quasilocal
fields to the study of Brauer groups do not restrict to stable
fields with Henselian valuations (see Section 6) and to the
corestriction mapping. In what follows, we use the abbreviation
tr-degree of the transcendency degree of an arbitrary field
extension.
\par
\medskip
{\bf Proposition 4.1.} {\it Let $E$ and $M _{1}, \dots , M _{s}$ be
fields, and for $j = 1, \dots , s$, suppose that $M _{j} \in {\rm
Gal}(E)$, $[M _{j}\colon E] = n _{j}$, $N _{M _{j}/E} = N _{M
_{j}/E} (\widehat X _{j})$ is a norm form of $M _{j}/E$ in an $n
_{j}$-tuple $\widehat X _{j} = (X _{j,1}, \dots , X _{j,n _{j}})$ of algebraically independent variables over $E$, $c _{j} \in E ^{\ast
}$, and $\Lambda (M _{j}/E; c _{j})$ is the fraction field of the
quotient ring $E[\widehat X _{j}]/(N _{M/E} - c _{j})$. Then the
field $\Lambda = \widetilde \otimes _{j=1} ^{s} \Lambda (M _{j}/E; c _{j})$,
where $\widetilde \otimes = \widetilde \otimes _{E}$, has the following properties:}
\par
(i) $\Lambda $ {\it is regular over $E$ and $(M \otimes _{E} \Lambda
)/M$ is rational of {\rm tr}-degree $\sum _{j=1} ^{s} (n _{j} - 1)$;}
\par
(ii) Br$(\Lambda /E)$ {\it equals the group $B _{c} (M/E) = \langle
[(L _{j}/E,$ $\sigma _{L _{j}/E}, c _{j})]$, $j = 1, \dots ,
s\rangle $, where $L _{j}$ runs across the set of cyclic extensions
of $E$ in $M _{j}$, and $\sigma _{L _{j}/E}$ is a generator of $G(L
_{j}/E)$. In particular, if $c _{j} \in N(M _{j}/E)$, $j = 1, \dots
, s$, then {\rm Br}$(\Lambda /E) = \{0\}$.}
\par
\medskip
{\it Proof.} (i): For each index $j$, $M _{j} [\widehat X _{j}]/(N
_{M/E} - c _{j})$ is a domain, its fraction field $\Theta _{j}$ is
rational over $M _{j}$ of tr-degree $n _{j} - 1$, and $G(M _{j}/E)$
acts on $\Theta _{j}$ as a group of automorphisms. Therefore,
$\Theta _{j}/M _{j}$ is regular and $\Lambda (M _{j}/E; c _{j})$ is
the fixed field of $G(M _{j}/E)$ in $\Theta _{j}$. In view of [21,
Ch. VII, Proposition 20], this implies $\Theta _{j} \in {\rm
Gal}(\Lambda (M _{j}/E; c _{j}))$, $G(\Theta _{j}/\Lambda (M _{j}/E;
c _{j})) \cong G(M _{j}/E)$ and $M _{j} \cap \Lambda (M _{j}/E; c
_{j}) = E$. It is now clear that $\Lambda (M _{j}/E; c _{j})/E$ is
regular and $\Theta _{j} \cong _{E} M _{j} \otimes _{E} \Lambda (M
_{j}/E; c _{j})$, so Proposition 4.1 (i) can easily be proved by
induction on $s$.
\par
(ii) If $s \ge 2$, then $\Lambda $ is $E$-isomorphic to $\widetilde
\otimes _{j=1} ^{s-1} \Lambda _{j} ^{\prime }$, where $\widetilde
\otimes =$ $\widetilde \otimes _{\Lambda (M _{s}/E; c _{s})}$ and
$\Lambda _{j} ^{\prime } = \Lambda (M _{j}/E; c _{j}) \otimes _{E}
\Lambda (M _{s}/E; c _{s})$, $j = 1, \dots ,$ $s - 1$. Thus our
proof reduces to the case of $s = 1$.
\par
To simplify notation, put $M = M _{1}$, $c = c _{1}$, $n = n _{1}$
and $X _{u} = X _{1,u}$, for $u = 1, \dots , n$. The extensions of
$E$ considered in the rest of our proof are assumed to lie in
$I(\overline \Theta /E)$, for some algebraically closed extension
$\overline \Theta $ of $E$ of countable tr-degree. In what follows,
$B = \{\xi _{u}\colon \ u = 1, \dots , n\}$ denotes the basis of
$M/E$ with respect to which $N _{M/E}$ is defined, $Y = \sum _{j=1}
^{n} \xi _{j}X _{j}$, and $\bar f$ stands for the image in
$M\Lambda (M/E; c)$ of any polynomial $f \in M[X _{1}, \dots , X
_{n}]$. Recall that $M\Omega \in {\rm Gal}(\Omega )$ and $G(M\Omega
/\Omega ) \cong G(M/E)$, for any $\Omega \in I(\overline \Theta /E)$
with $\Omega \cap M = E$. In addition, then $I(M/E)$ is mapped
bijectively on $I(M\Omega /\Omega )$, by the rule $D \to D\Omega $.
Moreover, $D \in {\rm Gal}(E)$ if and only if $D\Omega \in {\rm
Gal}(\Omega )$; when this is the case, one may identify when
necessary $G(D/D _{0})$ with $G(D\Omega /D _{0}\Omega )$, for any $D
_{0} \in I(D/E)$. Note that the norm map $N _{D _{0}\Omega }
^{D\Omega }$ extends $N _{D _{0}} ^{D}$. Let $R \in I(M/E)$,
$[R\colon E] = r$, $\psi _{1}, \dots , \psi _{r}$ the $E$-embeddings
of $R$ into $M$, and for each index $k \le r$, let $\tilde \psi
_{k}$ be an automorphism of $M$ extending $\psi _{k}$, and $\rho
_{k} = \prod _{\sigma \in G(M/R)} (\sigma \tilde \psi _{k} ^{-1})
(\bar Y)$. It is easily seen that $\bar \rho _{1}, \dots , \bar \rho
_{r-1}$ are algebraically independent over $R$, $\prod _{k=1} ^{r}
\bar \rho _{k} = c$ and $\bar \rho _{k} \in R\Lambda (M/E; c)$, $k =
1, \dots r$. Hence, the above observations and the transitivity of
norms in towers of finite extensions lead to the conclusion that
$R\Lambda (M/E; c)$ is $R$-isomorphic to $\widetilde \otimes _{R}$
of the fields $\Lambda (\widetilde M/\widetilde R; c)$ and $\Lambda
(\widetilde M/\widetilde R; \bar \rho _{k})$, $k = 1, \dots , r -
1$, where
$\widetilde R = R(\bar \rho _{1}, \dots , \bar \rho _{r-1})$ and
$\widetilde M = M\widetilde R$. This enables us to complement
Proposition 4.1 (ii) as follows:
\par
\medskip
(4.1) The sum of groups Br$(\Lambda (\widetilde M/\widetilde R;
c)/\widetilde R)$ and Br$(\Lambda (\widetilde M/\widetilde R; \rho
_{k})/\widetilde R)$, $k = 1, \dots ,$ $r - 1$, is direct and equal
to Br$(R\Lambda (M/E; c)/\widetilde R)$. Also, Br$(R\Lambda
(M/E; c)/R) =$ ${\rm Br}(\Lambda (M/R; c)/R)$.
\par
\medskip
Note that $c \in N(M\Lambda (M/E; c)/\Lambda (M/E;  c))$, since $N
_{\Lambda (M/E; c)} ^{M\Lambda (M/E; c)}$ extends $N _{E} ^{M}$. In
view of [27, Sect. 15.1, Proposition b], this implies that $B _{c}
(M/E) \subseteq {\rm Br}((M\Lambda (M/E; c))/$ $\Lambda (M/E; c))$.
We prove the converse implication. By [14, Theorem 1], there is
$B(M/E) \in I(\overline \Theta /E)$ which is $E$-isomorphic to
$\widetilde \otimes _{E}$ of function fields of Brauer-Severi
$E$-varieties, such that Br$(B(M/E)/E)$ $= B _{c} (M/E)$. This
ensures that $E$ is algebraically closed in $B(M/E)$. Therefore,
considering $MB(M/E)/B(M/E)$ instead of $M/E$, one obtains that
Proposition 4.1 (ii) will follow, if we show that Br$((M\Lambda
(M/E;$ $c))/\Lambda (M/E; c)) = \{0\}$ in case $B _{c} (M/E) =
\{0\}$. Applying Theorem 1.2 (i) and the concluding statement of
Theorem 1.2 (ii) as in the proof of Proposition 2.6 in [10] (with
$E$ instead of $E _{0}$), and using Theorem 1.1, one sees further
that it suffices to prove the final assertion of Proposition 4.1
(ii). Suppose first that $[M\colon E] = p$, for some $p \in
\overline P$. We show that Br$(\Lambda (M/E; c)/E) = \{0\}$ by
proving the following statement:
\par
\medskip
(4.2) If $c \in N(M/E)$, then $\Lambda (M/E; c)/E$ is rational of
tr-degree $p - 1$.
\par
\medskip
It is sufficient to establish (4.2) in the case of $c = 1$. Fix a
generator $\sigma $ of $G(M/E)$ and put $Y _{i} = \sum _{i=1} ^{p}
\sigma ^{i-1} (\xi _{i})X _{i}$. As $N _{\Lambda (M/E; 1)}
^{M\Lambda (M/E; 1)} (\overline Y _{1}) = 1$, Hilbert's Theorem 90
yields $\overline Y _{1} = \overline Z\sigma (\overline Z) ^{-1}$,
for some $Z \in M[X _{1}, \dots , X _{p}]$. Clearly, $\lambda
\overline Z\sigma (\lambda \overline Z) ^{-1} = \overline Y _{1}$,
for each $\lambda \in \Lambda (M/E; c) ^{\ast }$. Denote by $W _{1}
(\lambda ), \dots ,$ $W _{p} (\lambda )$ the coordinates of $\lambda
\overline Z$ in $\Lambda (M/E; 1)$ with respect to $B$.
Observe that $\lambda $ can be fixed so that $W _{p} (\lambda ) \in
E$. Since $\overline Y _{1}, \dots , \overline Y _{p-1}$ form a
tr-basis of $M\Lambda (M/E; c)/M$, the choice of $\lambda
$ guarantees that $W _{1} (\lambda ), \dots , W _{p-1} (\lambda )$
are algebraically independent over $E$. Note finally that the
equality $\overline Y _{1} = \overline Z\sigma (\overline Z) ^{-1}$
implies that $\bar X _{i} \in E(W _{1} (\lambda ), \dots , W
_{p-1} (\lambda ))$, for $i = 1, \dots p$, so it follows that
$\Lambda (M/E; c) = E(W _{1} (\lambda ), \dots , W _{p-1} (\lambda
))$, proving (4.2).
\par
Assume now that $[M\colon E] = p ^{k}$, for some $p \in P(E)$ and $k
\in \hbox{\Bbb N}$, $k \ge 2$. Proceeding by induction on $k$ and
using the fact that finite extensions of $E$ are quasilocal, one
obtains that it suffices to prove Proposition 4.1 (ii) when $c \in
N(M/E)$, under the hypothesis that its conclusion holds in general,
for every Galois extension $\Phi /\Phi _{0}$ of degree $p ^{u} < p
^{k}$, where $u \in \hbox{\Bbb N}$. This ensures the validity of
(4.1) for $M _{0} ^{\prime }/R _{0}$ whenever $M _{0} ^{\prime } \in
{\rm Gal}(E)$, $E \subseteq R _{0} \subseteq M _{0} ^{\prime }
\subseteq M$ and $[M _{0} ^{\prime }\colon R _{0}] < p ^{k}$. It is
well-known (see e.g. [21, Ch. I, Sect. 6; Ch. VIII]) that Gal$(E)
\cap I(M/E)$ contains fields $R _{1}, \dots , R _{k}$, such that $[R
_{j}\colon E] = p ^{j}$, $j =  1, \dots , k$, and $R _{j-1} \subset
R _{j}$, for $j \ge 2$. Put $H _{j} = G(M/R _{j})$, $y _{j} = (\prod
_{h _{j} \in H _{j}}h _{j} (\overline Y))$, $S _{j} = \{\tau _{j} (y
_{j})\colon \ \tau _{j} \in G(R _{j}/E)\}$ and $\widetilde R _{j} =
R _{j} (S _{j})$, for $j = 1, \dots , k$. It is easy to see that
$G(R _{j}/E)$ is an automorphism group of $\widetilde R _{j}$ whose
fixed field, say $N _{j}$, is $E$-isomorphic to $\Lambda (R _{j}/E;
c)$ and satisfies the equality $N _{j} (S _{j}) = \widetilde R _{j}$
($S _{j}$ will be viewed as a standard generating set of $\widetilde
R _{j}/R _{j}$ in $\Lambda (M/E; c)$). To prove that Br$(\Lambda
(M/E; c)/E) = \{0\}$ we show (setting $B _{0} = \Phi _{0} = E$) the
existence of a tower of field extensions $T _{1}, B _{1}, \Phi _{1},
\dots , T _{k-1}, B _{k-1}, \Phi _{k-1}, T _{k}$ of $N _{1}$
satisfying the following conditions, for each index $j \ge 1$:
\par
\medskip
(4.3) (i) $T _{j}/\Phi _{j-1}$ is rational of tr-degree $p ^{j}$
and contains $N _{j}$ as a subfield; more precisely, $T _{j}$ is a
transfer over $R _{j}\Phi _{j-1}$ to $\Phi _{j-1}$ of a rational
function field $T _{j} ^{\prime } = \Phi _{j-1} (Z _{j,1})$, such
that the compositum of $R _{j}$ and the transfer $T _{j,(j-1)}
^{\prime }$ of $T _{j} ^{\prime }$ over $(R _{j}\Phi _{j-1})/(R
_{j-1}\Phi _{j-1})$ includes the set $S _{j,(j-1)} = \{\tau _{j}
^{\prime } (y _{j})$, $\tau _{j} ^{\prime } \in G(R _{j}/R _{j-1}),
\tau _{j} ^{\prime } \neq 1\}$. Specifically, the union $\{Z
_{j,1}\} \cup S _{j,(j-1)}$ is a tr-basis of $(R _{j}T _{j,(j-1)}
^{\prime })/(R _{j}\Phi _{j-1})$.
\par
(ii) $B _{j}$ is an $(R _{j}T _{j})/T _{j}$-transfer of an
extension $B _{j} ^{\prime }$ of $R _{j}T _{j}$ isomorphic to a
tensor compositum over $R _{j}T _{j}$ of function fields of
Brauer-Severi $(R _{j}T _{j})$-varieties, and such that Br$(B _{j}
^{\prime }/R _{j}T _{j})$ equals the sum $\Gamma _{j}$ of the
images of the groups Br$(\Lambda (M\widetilde R _{j}/\widetilde R
_{j}; s _{j}))$, $s _{j} \in S _{j}$, under $\rho _{\widetilde R
_{j}/(R _{j}T _{j})}$. In particular, $T _{j}$ is algebraically
closed in $B _{j}$.
\par
(iii) Br$(B _{j}/T _{j})$ equals the image of $\Gamma _{j}$
under Cor$_{(R _{j}T _{j})/T _{j}}$, and $\Gamma _{j} \cap {\rm
Br}(\Phi _{j-1}) = \{0\}$.
\par
(iv) $\Phi _{j}$ is $p$-quasilocal, $B _{j} \subseteq \Phi _{j}$, $B
_{j}$ is algebraically closed in $\Phi _{j}$ and Br$(\Phi _{j}/B
_{j})$ $= \{0\}$.
\par
\medskip
The first part of (4.3) (iii) and the second half of (4.3) (ii) are
implied by the first part of (4.3) (ii) and the following lemma.
\par
\medskip
{\bf Lemma 4.2.} {\it With assumptions being as in Proposition 2.8,
suppose that $M \in {\rm Gal}(F)$. Then $M\Lambda $ is $\widetilde
\otimes _{M}$ of function fields of Brauer-Severi $M$-varieties and
$\Lambda ^{\prime }$ embeds canonically in $M\Lambda $ as an
$M$-subalgebra. Moreover, $E$ is algebraically closed in $\Lambda $
and {\rm Br}$(\Lambda /E)$ coincides with the images of {\rm
Br}$(\Lambda ^{\prime }/M)$ and {\rm Br}$(M\Lambda /M)$ under {\rm Cor}$_{M/E}$.}
\par
\medskip
{\it Proof.} The assertion about Br$(\Lambda /E)$ has been deduced
from Lemma 2.4 in [10, Sect. 2]. The other conclusions of the lemma
follow from its assumptions and Galois theory (see the beginning of
[30, Sect. 3], for more details).
\par
\medskip
Statement (4.3) (ii) and the properties of cyclic algebras described
in [27, Sect. 15.1] guarantee that $S _{j} \subset N(I _{j}B _{j}/R
_{j}B _{j})$, for every $I _{j} \in I(M/R _{j})$ cyclic over $R
_{j}$. As $B _{j}$ is separably closed in $\Phi _{j}$, this means
that $S _{j} \subset N(I _{j}\Phi _{j}/R _{j}\Phi _{j})$, so it
follows from Galois theory and the $p$-quasilocal property of $\Phi
_{j}$ (apply (2.3) (ii)) that $S _{j} \subset N(M\Phi _{j}/R
_{j}\Phi _{j})$, for each admissible $j$. These facts enable one to
deduce (4.3) (i) from the following lemma.
\par
\medskip
{\bf Lemma 4.3.} {\it Assume that $E$ is a field, $R$, $L$ and $M$
lie in {\rm Gal}$(E)$, $R \subset L \subseteq M $, $[L\colon R] =
p$ and $[M\colon E] = p ^{k}$, for some $p \in P(E)$, $k \in
\hbox{\Bbb N}$. Suppose also that the inductive hypothesis holds
and fix an element $\rho \in N(M/R)$ so that $c = N _{E} ^{R} (\rho
)$. Then $\Lambda (L/R; \rho ) \in I(\Omega ^{\prime }/R)$, for a
rational extension $\Omega ^{\prime }$ of $R$ satisfying the
following conditions:}
\par
(i) $\Omega ^{\prime }$ {\it is rational over $\Lambda (L/R; \rho )$
and has {\rm tr}-degree $p$ over $R$; also, $\Omega ^{\prime }$ is
the $L/R$-transfer of a rational extension $\Omega _{1}$ of $L$ in
one indeterminate;}
\par
(ii) {\it The $L/E$-transfer $\Omega $ of $\Omega _{1}$ is rational
over $E$ of {\rm tr}-degree $[L\colon E]$;}
\par
(iii) {\it The $R$-algebra $\widetilde \otimes _{R} \Lambda (L/R;
\tau (\rho ))$, where $\tau $ runs across $G(R/E)$, is isomorphic to
a field $\Lambda _{R} \in I(R\Omega /R)$, such that $R\Omega
/\Lambda _{R}$ is rational of {\rm tr}-degree $[R\colon E]$.}
\par
\medskip
{\it Proof.} It is clearly sufficient to consider the special case
where $\rho = 1$. Identifying $G(L/R)$ with $G(L\Lambda (L/R;
1)/\Lambda (L/R;$ $1))$, fix a generator $\sigma $ of $G(L/R)$ and
take elements $y _{i} \in L\Lambda (L/R; 1)\colon \ i = 1, \dots ,
p$, so that $\prod _{i=1} ^{p} y _{i} = 1$, $\{y _{1}, \dots , y
_{p-1}\}$ is a tr-basis of $\Lambda (L/R; 1)/L$ and $y _{i} = \sigma
^{i-1} (y _{1})$, for each index $i$. Note also that $L\Lambda (L/R;
1)$ is a subfield of a rational function field $L ^{\prime } = L(z
_{1}, \dots , z _{p})$ with generators subject to the relations $z
_{i}/z _{i+1} ^{-1} = y _{i}$, $i = 1, \dots , p - 1$. This
guarantees the existence of an $E$-automorphism $\tilde \sigma $ of
$L ^{\prime }$ extending $\sigma $, and such that $\tilde \sigma (z
_{p}) = z _{1}$ and $\tilde \sigma (z _{i}) = z _{i+1}$, $i < p$.
Therefore, $G(L/R)$ can be viewed as a group of automorphisms of $L
^{\prime }$ whose fixed field $\Omega ^{\prime }$ includes $\Lambda
(L/R; 1)$. Since, by Hilbert's Theorem 90, $y _{p} = w _{p}\sigma (w
_{p}) ^{-1}$, for some $w _{p} \in L\Lambda (L/R; 1)$, it follows
that $z _{p} = zw _{p}$, for some $z \in \Omega ^{\prime }$
transcendental over $\Lambda (M _{0}/E; 1)$. Thus it becomes clear
that $\Omega ^{\prime } = \Lambda (L/R; 1) (z)$ and $L ^{\prime } =
L\Lambda (M _{0}/E; 1) (z)$. Putting $\Omega _{1} = L(z _{1})$ and
applying [30, Lemmas 3.1 and 3.2], one completes the proof of (i)
and (ii). Statement (iii) is implied by (i), (ii) and the definition
of the transfer map.
\par
\medskip
The latter assertion of (4.3) (iii) follows from the former one and
our next lemma.
\par
\medskip
{\bf Lemma 4.4.} {\it In the setting of Lemma 4.3, let $\Lambda 
_{R}$ be the field $\widetilde \otimes _{\tau \in G(R/E)} \Lambda 
(M/R; \tau (\rho ))$, where $\widetilde \otimes = \widetilde \otimes 
_{R}$. For each $\tau \in G(R/E)$, identify $\Lambda (L/R;$ $\tau 
(\rho ))$ with its canonical $R$-isomorphic copy in $\Lambda (M/R;$ 
$\tau (\rho ))$, and denote by $S _{\tau (\rho )}$ the standard 
generating set of $L\Lambda (L/R;$ $\tau (\rho ))/L$ in $\Lambda 
(M/R; \tau )$. Also, let $\Sigma _{L\Lambda _{R}}$ be the sum of the 
images of
\par \noindent
{\rm Br}$(\Lambda (M(S _{\tau (\rho )})/L(S _{\tau (\rho 
)});$ $s _{\tau (\rho )})$ under $\rho _{L(S _{\tau (\rho 
)})/L.\Lambda _{R}}$, where $\tau $ runs across $G(R/E)$ and, for 
each $\tau $, $s _{\tau (\rho )}$ runs across $S _{\tau (\rho )}$. 
Then the image $\Delta _{E}$ of $\Sigma _{L\Lambda _{R}}$ under {\rm Cor}$_{(L\Omega )/\Omega } \circ \rho _{L\Lambda _{R}/L\Omega }$ 
intersects trivially with {\rm Br}$(E)$.}
\par
\medskip
{\it Proof.} The assertion is obvious in the case where $E$ is
finite, since then Br$(E) = \{0\}$, by Wedderburn's theorems.
Suppose further that $E$ is infinite, denote by $\Delta _{L}$ the
image of $\Sigma _{L\Lambda _{R}}$ under $\rho _{(L\Lambda
_{R})/(L\Omega )}$, let $G(L/E) = \{\varphi _{1}, \dots , \varphi
_{l}\}$, $\varphi _{1} = 1$, $[L\colon E] = l$, and take a tr-basis
$Z _{1}, \dots , Z _{l}$ of $L\Omega /\Omega $ so that $L(Z _{1}) =
\Omega _{1}$ and $Z _{j} = \varphi _{j} (Z _{1})$, for each index
$j$. Statement (4.1) and the inductive hypothesis on $M/E$ imply
that $\Delta _{L} \subseteq {\rm Br}(M\Omega /L\Omega )$ and $\Delta
_{E} \subseteq {\rm Br}(M\Omega /\Omega )$. In addition, it follows
from Lemma 4.3 (ii) that Br$(E) \cap \Delta _{E} \subseteq {\rm
Br}(M/E)$. Denote by $A _{L} ^{l}$ the $l$-dimensional affine
$L$-space, and for each $d \in \hbox{\Bbb N}$, let $A _{L/E} ^{d} =
\{(a _{1} ^{d}, \dots , a _{l} ^{d})\colon \ a _{i} \in L ^{\ast } \
{\rm and} \ a _{i} = \varphi _{i} (a _{1}), i = 1, \dots , l\}$.
Observing that $L$ is algebraic over $E(Z _{1} ^{d}, \dots , Z _{l}
^{d})$, and using (4.1), the infinity of $E$ and the definition of
the corestriction mapping, one proves the following:
\par
\medskip
(4.4) (i) The sets $A _{E} ^{l}$ and $A _{L/E} ^{d}\colon \ d \in
\hbox{\Bbb N}$, are dense in $A _{L} ^{l}$, in the sense of Zariski;
\par
(ii) $A _{L/E} ^{l}$ possesses a subset $U \neq \phi $,
Zariski-open in $A _{L} ^{l}$ and such that each specialization of
$(Z _{1}, \dots , Z _{l})$ into $U$ induces group homomorphisms
$\pi _{L}\colon \ ({\rm Br}(M/L) + \Delta _{L}) \to {\rm Br}(L)$
and $\pi _{E}\colon \ ({\rm Br}(M/E) + \Delta _{E}) \to {\rm
Br}(E)$ satisfying the equality Cor$_{L/E} \circ \pi _{L} = \pi
_{E} \circ {\rm Cor}_{(L\Omega )/\Omega )}$ and acting as the
identity on Br$(M/L)$ and Br$(M/E)$, respectively.
\par
\medskip
Statement (4.4) (i) and [27, Sect. 15.1, Proposition b] ensure the
existence of many specializations for which $\Delta _{L} \subseteq
{\rm Ker}(\pi _{L})$. This implies that $\Delta _{E} \subseteq {\rm
Ker}(\pi _{E})$ and Br$(E) \cap \Delta _{E} = \{0\}$, so Lemma 4.4
is proved.
\par
\medskip
We are now in a position to complete the proof of (4.3) and
Proposition 4.1. The existence of $B _{j} ^{\prime }$ satisfying the
conditions in the first part of (4.3) (ii) is obtained by applying
[14, Theorem 1]. Statement (4.3) (iv) follows from Lemma 2.4, so
(4.3) holds and Br$(T _{k}/E) = \{0\}$. Since, by (4.3) (i), $N
_{k} \in I(T _{k}/E)$ and $N _{k}$ is $E$-isomorphic to $\Lambda
(M/E; c)$, this proves Proposition 4.1 (ii) in case $M \subseteq E
(p)$.
\par
For the proof in general, take $G _{p}$ and $E _{p}$ as in Lemma
3.6. By Lemma 2.2,
\par \noindent
$c \in N(M/E _{p})$, so (4.1) and the established special case of
our assertion imply that Br$(E _{p}\Lambda (M/E; c))/E _{p}) =
\{0\}$. It is now easy to see that Br$(\Lambda (M/E; c)/E) \cap {\rm
Br}(E) _{p} = \{0\}$, for every $p \in \overline P$ (see [27, Sect.
13.4]), which completes the proof of Proposition 4.1 (ii).
\par
\medskip
{\bf Remark 4.5.} Propositions 2.6 and 4.1 imply that if $E$ is a
field, $M \in {\rm Gal}(E)$ and $R \in I(M/E)$, then Cor$_{R/E}$
maps Br$(\Lambda (M/R;$ $c)/R)$ into Br$(\Lambda (M/E; c)/E)$, for
each $c \in E ^{\ast }$. Moreover, it follows from the RC-formula
that Br$(\Lambda (M/E; c)/E)$ equals the image of Br$(\Lambda (M/R;
c)/R)$ under Cor$_{R/E}$, provided that ${\rm g.c.d.}([R\colon E],
[M\colon R]) = 1$.
\par
\medskip
{\bf Proposition 4.6.} {\it Let $E$ be a field, $\Omega $ an
algebraically closed extension of $E _{\rm sep}$, $s$ a positive
integer, $\Sigma $, $\Lambda $ and $M _{1}, \Lambda _{1}, \dots M
_{s}, \Lambda _{s}$ lie in $I(\Omega /E)$ so that $\Lambda = \Lambda
_{1} \dots \Lambda _{s}$, $\Lambda \cap \Sigma _{\rm sep} = E$ and
$\Lambda \cong \widetilde \otimes _{j=1} ^{s} \Lambda _{j}$, where
$\widetilde \otimes = \widetilde \otimes _{E}$ and $\cong $ is an
$E$-isomorphism. For each $j$, assume that $M _{j} \in {\rm
Gal}(E)$, put $\Sigma _{j} = \Sigma \cap M$ and let $\Lambda _{j}/E$
be of one of the following types:}
\par
(i) $\Lambda _{j} = \Lambda (M _{j}/E; c _{j})${\it , for some $c
_{j} \in E ^{\ast }$;}
\par
(ii) $\Lambda _{j}/E$ {\it is an $M _{j}/E$ transfer of $\otimes _{M
_{j}}$ of function fields of Brauer-Severi $M _{j}$-varieties, such
that ${\rm Br}(\Lambda _{j}/M _{j})$ is a submodule of {\rm Br}$(M
_{j})$ over the integral group ring $\hbox{\Bbb Z} [G(M _{j}/E)]$.
\par \noindent
Suppose also that $\widetilde \Delta _{j} = {\rm Br}(\Lambda
_{j}\Sigma _{j}/\Sigma _{j})$ in case (i), and let $\widetilde
\Delta _{j}$ be the image of {\rm Br}$(\Lambda _{j}M _{j}/M _{j})$
under {\rm Cor}$_{M _{j}/\Sigma _{j}}$, otherwise. Then $E$ and
$\Sigma $ are algebraically closed in $\Lambda $ and $\Lambda \Sigma
$, respectively, and {\rm Br}$(\Lambda \Sigma /\Sigma )$ equals the
sum of the images $\Delta _{j}$ of $\widetilde \Delta _{j}$ under
$\rho _{\Sigma _{j}/\Sigma }$, for $j = 1, \dots s$.}
\par
\medskip
{\it Proof.} It is clear from Galois theory and the condition
$\Lambda \cap \Sigma _{\rm sep} = E$ that if $F \in I(\Sigma /E)$,
then $F$ is algebraically closed in $\Lambda F$. For the rest of our
proof, suppose first that $s = 1$. It follows from Galois theory and
the definitions of $\Lambda $ and of the transfer map that $\Lambda
/E$ and $\Lambda \Sigma /\Sigma $ are of one and the same type
relative to $M _{1}/E$ and $M _{1}\Sigma /\Sigma $, respectively. In
addition, if $\Lambda /E$ is of type (ii), then direct calculations
show that the mappings of Br$(M _{1})$ into Br$(\Sigma )$ defined by
the rules Cor$_{M _{1}\Sigma /\Sigma } \circ \rho _{M _{1}/M
_{1}\Sigma }$ and $\rho _{\Sigma _{1}/\Sigma } \circ {\rm Cor}_{M
_{1}/\Sigma _{1}}$ coincide. Applying now Proposition 4.1 and (4.1)
in case $\Lambda /E$ is of type (i), and using Lemma 4.2, otherwise,
one proves that Br$(\Lambda \Sigma /\Sigma ) = \Delta _{1}$, as
claimed. It remains to be seen that the concluding assertion of
Proposition 4.6 holds in the case of $s \ge 2$. Consider the fields
$\Lambda $, $\widetilde E = \Lambda _{1} \dots $ $\Lambda _{s-1}$,
$\widetilde \Sigma = \Sigma \widetilde E$ and $M _{s}\widetilde E$
instead of $\Lambda $, $E$, $\Sigma $ and $M _{1}, \dots , M _{s}$,
respectively. It is not difficult to deduce from Galois theory and
well-known properties of tensor products and of function fields of
Brauer-Severi varieties (see [27, Sects. 9.2 and 9.4] and the proof
of [14, Theorem 1]) that Br$(\Lambda \widetilde \Sigma /\widetilde
\Sigma )$ is determined in accordance with Proposition 4.1, for $s =
1$, and equals the image of Br$(\Lambda _{s}\Sigma _{s}/\Sigma
_{s})$ under $\rho _{\Sigma _{s}/\widetilde \Sigma }$. Proceeding
now by induction on $s$ and applying the inductive hypothesis to
Br$(\widetilde \Sigma /\Sigma )$, we complete our proof.
\par
\vskip0.6truecm
\centerline{\bf 5. Proof of Theorem 1.2}
\par
\medskip
Let $E _{0}$ be an arbitrary field. Theorem 1.2 will be proved by
constructing $E$ as a union of a certain tower of fields $E
_{n}\colon \ n \in \hbox{\Bbb N}$, such that $E _{0} \subset E _{1}$
and $E _{n-1}$ is algebraically closed in $E _{n}$, for every index
$n$. It should be emphasized that the proof of Theorem 1.2 (i) and
of the latter part of Theorem 1.2 (ii) in the special case where
$\chi = {\rm Fin}$ does not use Proposition 2.8 in full generality
and is independent of Propositions 2.6, 4.1 and 4.6 (i) (but relies
on Lemma 4.2 and Propositions 3.4 (c), 4.6 (ii)). In order to ensure
generally that our construction has the desired properties we also
need the following lemmas.
\par
\medskip
{\bf Lemma 5.1.} {\it Let $E$ and $M$ be fields, and $\chi $ an
abelian closed class, such that ${\rm Nil} \subseteq \chi $, $M \in
{\rm Gal}(E)$, $G(M/E) \not\in \chi $ and $G(L/E) \in \chi $, for
$L \in {\rm Gal}(E) \cap I(M/E)$, $L \neq M$. Then:}
\par
(i) $G(M/E)$ {\it is simple or has a unique minimal normal subgroup
$G _{0}$; in the former case, $G(M/E) \not\in {\rm Sol}$;}
\par
(ii) {\it If $G(M/E)$ is not simple, then $G _{0} \in {\rm Ab}$ if
and only if $G(M/E) \in {\rm Sol}$; in this case, $\chi = {\rm
Nil}$;}
\par
(iii) {\it If $G(M/E) \in {\rm Sol}$, then $G _{0} \subseteq
[G(M/E), G(M/E)]$ and $G _{0} \in {\rm Syl}_{p} (M/E)$, for some $p
\in \overline P$ not dividing the order of $G(M/E)/G _{0}$.}
\par
\medskip
{\it Proof.} Suppose for a moment that $G(M/E)$ has normal proper
subgroups $H _{1}$ and $H _{2}$, such that $H _{1} \cap H _{2} =
\{1\}$. Then Galois theory and our assumptions ensure that $G(M/E)/H
_{j} \in \chi $, for $j = 1, 2$. Hence, by the choice of $\chi $, it
contains $G(M/E)/H _{1} \times G(M/E)/H _{2}$, and since $G(M/E)$
embeds canonically into $G(M/E)/H _{1} \times G(M/E)/H _{2}$, this
requires that $G(M/E) \in \chi $, a contradiction proving Lemma 5.1
(i). In the rest of the proof, we may assume that $G(M/E)$ has a
unique minimal normal subgroup $G _{0}$. It is well-known that if
$G(M/E) \in {\rm Sol}$, then $G _{0} \in {\rm Ab}$ and $G _{0}$ is
of exponent $p \in \overline P$. Conversely, if $G _{0} \in {\rm
Ab}$, then $\chi $ could not be abelian closed (by Galois theory and
the assumptions on $M/E$), so the conditions on $\chi $ guarantee
that it equals ${\rm Nil}$. Since Sol is closed under the formation
of group extensions, this proves Lemma 5.1 (ii). Suppose finally
that $G(M/E) \in {\rm Sol}$ and fix a group $G _{p} \in {\rm
Syl}_{p} (M/E)$. As $G(M/E) \not\in {\rm Nil}$ and Nil is a
saturated group formation (in the sense of [34]), $G _{0}$ is not
included in the Frattini subgroup of $G(M/E)$. Hence, $G(M/E) = G
_{0}H$, for some maximal subgroup $H$ of $G(M/E)$. In view of Lemma
5.1 (ii), this means that $G _{0} \cap H = \{1\}$ and $H \cong
G(M/E)/G _{0}$. In particular, $H \in {\rm Nil}$, which implies that
$G _{p}$ is normal in $G(M/E)$. The centre $Z(G _{p})$ of $G _{p}$
is characteristic in $G _{p}$, so the obtained result shows that
$Z(G _{p})$ is normal in $G(M/E)$. As $G _{0} \cap Z(G _{p}) \neq
\{1\}$ (see [21, Ch. I, Sect. 6]), the minimality of $G _{0}$
implies that $G _{0} \subseteq Z(G _{p})$. It is now easily seen
that the group $H _{p} = G _{p} \cap H$ is normal in $G(M/E)$. Since
$H _{p} \cap G _{0} = \{1\}$, Lemma 5.1 (i) yields $H _{p} = \{1\}$.
Note finally that if $G _{0} \not\subseteq [G(M/E), G(M/E)]$, then
$G _{0}$ must be a direct summand in $G(M/E)$. This, however, means
that $G(M/E) \cong G _{0} \times G(M/E)/G _{0}$, which contradicts
the assumption that $G(M/E) \not\in {\rm Nil}$, and so proves Lemma
5.1 (iii).
\par
\medskip
{\bf Lemma 5.2.} {\it Let $\Phi $, $L$ and $M$ be fields, such that
$L, M \in {\rm Gal}(\Phi )$, and let $p \in \overline P$ be a
divisor of $[M\colon L \cap M]$. Suppose that $R _{p}$ and $Y
_{p}$ are the fixed fields of some groups $\widetilde H _{p} \in
{\rm Syl}_{p} (L/L \cap M)$ and $H _{p} \in {\rm Syl}_{p} (M/L \cap 
M)$, respectively, and $\delta _{p} \neq 0$ be an element of {\rm 
Ker}$(Y _{p}/L \cap M)$. Then $\rho _{Y _{p}/(R _{p}Y _{p})}
(\delta _{p}) \not\in {\rm Im}(R _{p}/R _{p}Y _{p})$.}
\par
\medskip
{\it Proof.} It is easily verified that $[R _{p}Y _{p}\colon L \cap 
M] = [R _{p}\colon L \cap M][Y _{p}\colon L \cap M]$. This implies 
that $p \not\vert [R _{p}Y _{p}\colon Y _{p}]$, whence $\rho _{Y
_{p}/(R _{p}Y _{p})} (\delta _{p}) \colon = \bar \delta _{p} \neq
0$. Also, it follows that Cor$_{(R _{p}Y _{p})/R _{p}}
(\bar \delta _{p})$ $= \rho _{(L \cap M)/R _{p}} ({\rm Cor}_{E
_{p}/(L \cap M)} (\delta _{p})) = 0$. On the other hand, if $\bar
\delta _{p} = \rho _{R _{p}/(R _{p}Y _{p})} (\tilde \delta _{p})$,
for some $\tilde \delta _{p} \in {\rm Br}(R _{p})$, then the
RC-formula yields Cor$_{(R _{p}Y _{p})/R _{p}} (\bar \delta _{p}) =
\tilde \delta _{p} ^{m _{p}} \neq 0$, where $m _{p} = [Y _{p}\colon
L \cap M]$. The obtained contradiction proves our assertion.
\par
\medskip
Let $E$ be a field and $M \in {\rm Gal}(E)$. Before stating our next
lemma, we denote by $N(M/E) _{\rm cyc}$ the intersection of the norm
groups of cyclic extensions of $E$ in $M$. By Proposition 4.1, $N(M/E)
_{\rm cyc} = \{c \in E ^{\ast }\colon \ {\rm Br}(\Lambda (M/E; c)/E)
= \{0\}\}$.
\par
\medskip
{\bf Lemma 5.3.} {\it Assume that $E$, $M$ and $\chi $ satisfy the
conditions of Lemma 5.1, put $M _{0} = M \cap E _{\chi }$, take a
divisor $p \in \overline P$ of $[M\colon M _{0}]$, and suppose that
$H _{p} \in {\rm Syl}_{p} (M/M _{0})$, $Y _{p}$ is the fixed field
of $H _{p}$, and $c$ is an element of $E ^{\ast } \setminus N(M/Y
_{p}) _{\rm cyc}$. Let also $M _{1}, \dots , M _{s}$ be fields lying
in {\rm Gal}$(E)$, for some $s \in \hbox{\Bbb N}$, and let $\Omega $,
$\Lambda _{1}, \dots , \Lambda _{s}$ and $\Lambda $ be extensions of
$E$ associated with $M _{1}, \dots , M _{s}$ as in Proposition 4.6.
Then $c \not\in N(M\Lambda /Y _{p}\Lambda ) _{\rm cyc}$ in the
following cases:}
\par
(i) $(M _{1} \dots M _{s}) \cap M \subseteq M _{0}$ {\it and $\chi
\neq {\rm Nil}$;}
\par
(ii) $M _{i} \subseteq E _{\rm Nil}$ {\it and $\Lambda _{i} = \Lambda
(M _{i}/E; c _{i})$, for some $c _{i} \in E ^{\ast }$, and $i = 1,
\dots , s$;}
\par
(iii) {\it For each index $j$, $\Lambda _{j}$ is of type (ii) (in
the sense of Proposition 4.6) and
\par \noindent
{\rm Br}$(M _{j}\Lambda _{j}/M _{j}) \subseteq {\rm Ker}_{M _{j}/L
_{j}}$, where $L _{j} \in {\rm Gal}(E)$ is chosen so that $L _{j}
\subseteq M _{j} \subseteq L _{j} (p _{j})$, for some $p _{j} \in
\overline P$. In this case, {\rm Br}$(Y _{p}\Lambda /Y _{p}) =
\{0\}$.}
\par
\medskip
{\it Proof.} Arguing as in the concluding part of the proof of
Proposition 4.6, one reduces our considerations to the case of $s =
1$. As $E$, $M$ and $\chi $ satisfy the conditions of Lemma 5.1, we
have $M \subseteq M _{1}$ or $M \cap M _{1} \subseteq M _{0}$. Our
first objective is to prove Lemma 5.3 (iii). Observe that if $M
\subseteq M _{1} \subseteq L _{1} (p _{j})$, then Br$(Y _{p}\Lambda
/Y _{p}) = \{0\}$. Indeed, it follows from Galois theory and the
assumptions on $M/E$ that $M \subseteq L _{1}$ except, possibly, in
the case of $\chi = {\rm Nil}$ and $p _{1} = p$. When $\chi = {\rm
Nil}$ and $p _{1} = p$, we have $Y _{p} = M _{0} \subseteq L _{1}$.
It is therefore clear from Lemma 4.2 that Br$(Y _{p}\Lambda /Y _{p})
= \{0\}$, as required. The same assertion is implied in the case of
$L _{1} \cap M \neq M$ by Proposition 4.6, since then Cor$_{L
_{1}/(L _{1} \cap M)}$ is injective and $L _{1} \cap M \subseteq M
_{0} \subset Y _{p}$. Hence, by Galois theory and [27, Sect. 15.1,
Proposition b], $c \not\in N(M\Lambda /Y _{p}\Lambda ) _{\rm cyc}$,
so Lemma 5.3 (iii) is proved. Assume now that $M \cap M _{1}
\subseteq M _{0}$, $\widetilde H _{p} \in {\rm Syl}_{p} (M _{0}M
_{1}/M _{0})$ and $R _{p}$ is the fixed field of $\widetilde H
_{p}$. It follows from (4.1) and Proposition 4.1 that if $\chi \neq
{\rm Nil}$, then Br$(R _{p}Y _{p}\Lambda /R _{p}Y _{p})$ is the
image of Br$(R _{p}\Lambda /R _{p})$ under $\rho _{R _{p}/(R _{p}Y
_{p})}$. When $\chi = {\rm Nil}$, $M _{1} \dots M _{s} \subseteq M
_{0}$, which enables one to obtain similarly that Br$(K _{p}\Lambda
/K _{p}) \cap {\rm Br}(K _{p}) _{p}$ is the image of Br$(\Lambda /E)
\cap {\rm Br}(E) _{p}$ under $\rho _{E/K _{p}}$, for $K _{p} = R
_{p}, R _{p}Y _{p}$. On the other hand, Lemma 5.1 and Galois theory
show that Lemma 3.6 applies to $(M/E, p)$, if $\chi = {\rm Nil}$,
and to $(M/M _{0}, p)$, otherwise. In view of Lemma 5.2, these
observations prove Lemma 5.3.
\par
\medskip
Let now $E _{0}$ be an arbitrary field, $T _{0}$ a subgroup of Br$(E
_{0})$ and $R _{\rm Fin}$ a system of representatives of the
isomorphism classes of finite groups. Replacing, if necessary, $E
_{0}$ by its rational extensions of sufficiently large tr-degree,
one easily reduces (e.g., from [27, Sect. 19.6]) our considerations
to the special case where $T$ is a divisible hull of $T _{0}$. Note
further that $E _{0}$ has a regular PQL-extension $E _{0} ^{\prime
}$, such that Br$(E _{0} ^{\prime }/E _{0}) = \{0\}$ and Br$(E _{0}
^{\prime })$ is divisible (apply [10, Lemma 1.4], proved on the
basis of (2.2)), so one may assume for the proof that $T _{0} = T$.
It is known [39] that each profinite group $G$ is (continuously)
isomorphic to $G(L(G)/K(G))$, for some rational extension $L(G)/E$
of countable tr-degree and a suitably chosen $K(G) \in I(L(G)/E)$.
Since Br$(L(G)/E) = {\rm Br}(K(G)/E) = \{0\}$, this applied to the
case in which $G$ is a topological product of the groups in $R _{\rm
Fin}$, allows us to assume for the proof of Theorem 1.2 that all $G
\in {\rm Fin}$ are realizable as Galois groups over $E _{0}$. It
follows from the choice of $E _{0}$ that it has Galois extensions
$\Sigma _{0,1}$ and $\Sigma _{0,2}$ in $E _{0,{\rm sep}}$, such that
$\Sigma _{0,1} \cap \Sigma _{0,2} = E$ and each $G \in {\rm Fin}$ is
isomorphic to $G(Y _{j}/E)$, where $Y _{j} \in {\rm Gal}(E)$ and $Y
_{j} \subseteq \Sigma _{0,j}$, for $j = 1, 2$. Our objective is to
prove the existence of a quasilocal extension $E/E _{0}$ with the
properties required by Theorem 1.2. The field $E$ will be obtained
as a union $\cup _{n=1} ^{\infty } E _{n}$ of an inductively defined
tower of regular extensions of $E _{0}$. Suppose that the field $E
_{k}$ has already been defined, for some integer $k \ge 0$, and
denote by $T _{k}$ the image of $T _{0} = T$ under $\rho _{E _{0}/E
_{k}}$. As $T _{k}$ is divisible, it is a direct summand in Br$(E
_{k})$, i.e. Br$(E _{k})$ possesses a subgroup $T _{k} ^{\prime }$,
such that $T _{k} ^{\prime } + T _{k} = {\rm Br}(E _{k})$ and $T
_{k} ^{\prime } \cap T _{k} = \{0\}$. Hence, by [14, Theorem 1],
there is an extension $\Lambda _{k}$ of $E _{k}$, such that
Br$(\Lambda _{k}/E _{k}) = T _{k} ^{\prime }$ and $\Lambda _{k}$ is
presentable as $\widetilde \otimes _{E _{k}}$ of function fields of
Brauer-Severi $E _{k}$-varieties; in particular $\Lambda _{k}/E
_{k}$ is regular. Identifying $E _{k,{\rm sep}}$ with its $E
_{k}$-isomorphic copy in $\Lambda _{k,{\rm sep}}$, put $\Sigma
_{k,j} = \Sigma _{0,j}\Lambda _{k}$, for $j = 1, 2$. The regularity
of $\Lambda _{k}/E _{k}$ ensures that $E _{k}$ is algebraically
closed in $\Lambda _{k}$, so it follows from Galois theory that
$\Sigma _{k,1}/\Lambda _{k}$ and $\Sigma _{k,2}/\Lambda _{k}$ have
the same properties as $\Sigma _{0,1}/E _{0}$ and $\Sigma _{0,2}/E
_{0}$. Denote by $Z _{k}$ the extension of $\Lambda _{k}$ defined as
follows:
\par
\medskip
(5.1) ($\alpha $) If $\chi = {\rm Nil}$, then $Z _{k} = \Sigma
_{k,1} \cap \Lambda _{k,\chi '}$ in case $\chi ^{\prime } \neq {\rm
Nil}$, and $Z _{k} = \Lambda _{k}$, otherwise.
\par
($\beta $) If $\chi \neq {\rm Nil}$, then $Z _{k} = \Lambda
_{k,\chi }(\Sigma _{k,1} \cap \Lambda _{k,\chi '})$ in case $\chi
^{\prime } \neq \chi $, and $Z _{k} = \Lambda _{k,\chi }$, otherwise.
\par
\medskip
Let now $M _{k} \in {\rm Gal}(\Lambda _{k})$ and $W(M _{k})$ be a
tensor compositum over $M _{k}$ of function fields of Brauer-Severi
$M _{k}$-varieties, such that Br$(W(M _{k})/$ $M _{k}) = {\rm Ker}(M _{k}/\Lambda _{k})$. Denote by $W(M _{k}/\Lambda _{k})$ the $M
_{k}/\Lambda _{k}$-transfer of $W(M _{k})$ and put $W(M _{k})
^{\prime } = M _{k}W(M _{k}/\Lambda _{k})$. It is known that
Br$(W(M _{k}) ^{\prime }/$ $M _{k}) = {\rm Br}(W(M _{k})/M _{k})$,
and by Lemma 4.2, Br$(W(M _{k}/\Lambda _{k})/\Lambda _{k}) = \{0\}$.
Denote by $W _{k}$ the tensor compositum over $\Lambda _{k}$ of the
fields $W(M _{k}/\Lambda _{k})$, taken over all $M _{k} \in {\rm
Gal}(\Lambda _{k})$, $M _{k} \subseteq Z _{k}$. It is easily
obtained from Galois theory and case (ii) of Proposition 4.6 that
$\Lambda _{k}$ is algebraically closed in $W _{k}$ and Br$(W
_{k}/\Lambda _{k}) = \{0\}$. For each $p \in \overline P$, let $Q(W
_{k}) _{p}$ be the set of all pairs $\widetilde W _{k,p} = (W _{k}
^{\prime }, W _{k,p} ^{\prime }) \in {\rm Gal}(W _{k}) \times {\rm
Gal}(W _{k})$, for which $W _{k,p} ^{\prime } \in I(W _{k} ^{\prime
}/W _{k} ^{\prime })$. Replacing Ker$({\rm Cor}_{M _{k}/\Lambda
_{k}})$ by Ker$({\rm Cor}_{W _{k,p}'/W _{k}'}) \cap {\rm Br}(W
_{k,p} ^{\prime }) _{p}$, attach to each $\widetilde W _{k,p} \in
Q(W _{k}) _{p}$ a field extension $W(\widetilde W _{k,p})/W _{k,p}
^{\prime }$ in the same way as $W(M _{k})/M _{k}$ is associated with
$(\Lambda _{k}, M _{k})$, and let $W(\widetilde W _{k,p};$ $W _{k,p}
^{\prime }/W _{k})$ be the $W _{k,p} ^{\prime }/W _{k}$-transfer of $W(\widetilde W _{k,p})$. Consider the tensor compositum $\Theta
_{k}$ over $W _{k}$ of the fields $W(\widetilde W _{k,p};$ $W _{k,p}
^{\prime }/W _{k})\colon $ $p \in \overline P$, $\widetilde W _{k,p}
\in Q(W _{k}) _{p}$. By Proposition 4.6 (ii), $W _{k}$ is
algebraically closed in $\Theta _{k}$ and Br$(\Theta _{k}/W _{k}) =
\{0\}$. Observing further that ${\rm Sol}$ equals the intersection of
abelian closed group classes, and applying Lemmas 2.2 and 5.3, one
obtains the following result:
\par
\medskip
(5.2) If $E _{k} ^{\prime }/E _{k}$ and $\widetilde \chi $ satisfy
the conditions of Lemma 5.1, then $J _{k} \cap N(E _{k} ^{\prime
}/(E _{k} ^{\prime } \cap E _{k,\tilde \chi }))$ coincides with $J
_{k} \cap N(E _{k} ^{\prime }\Theta _{k}/(E _{k} ^{\prime } \cap E
_{k,\tilde \chi })\Theta _{k})$ in the following cases:
\par
(i) $\widetilde \chi = \chi ^{\prime } \cup {\rm Sol}$, $G(E _{k}
^{\prime }/(E _{k} ^{\prime } \cap E _{k,\chi '})) \not\in {\rm
Sol}$ and $J _{k} = E _{k} ^{\prime } \cap E _{k,\chi '}$.
\par
(ii) $\widetilde \chi = \chi \cup {\rm Sol}$, $G(E _{k} ^{\prime }/E
_{k}) \in \chi ^{\prime }$, $J _{k} = E _{k} ^{\prime } \cap E
_{k,\chi }$, $G(E _{k} ^{\prime }/J _{k}) \not\in {\rm Sol}$ and $E
_{k} ^{\prime } \subset \Sigma _{0,2}E _{k}$.
\par
(iii) $\widetilde \chi = \chi = {\rm Nil}$, $G(E _{k} ^{\prime }/E
_{k}) \in {\rm Sol}$, $J _{k} = E _{k}$, and in the case of $\chi
^{\prime } \neq \chi $, $E _{k} ^{\prime } \subset \Sigma _{0,2}E
_{k}$.
\par
\medskip
Putting $\widetilde \Theta _{k} = \Theta _{k}$ and $E _{k+1} =
\Theta _{k}$ in the cases of $\chi ^{\prime } = {\rm Fin}$ and $\chi
= {\rm Fin}$, respectively, we continue the presentation of the
inductive step in our construction with the definition of the
extension $\widetilde \Theta _{k}/\Theta _{k}$ in the case where
$\chi ^{\prime } \neq {\rm Fin}$. Let Un$(\Theta _{k}) = \{\Theta
_{k} ^{\prime } \in {\rm Gal}(\Theta _{k})\colon \ G(\Theta _{k}
^{\prime }/(\Theta _{k} ^{\prime } \cap \Theta _{k,\chi '})) \not\in
{\rm Sol}\}$, and let $\Omega _{k}/\Theta _{k}$ be a rational
extension with a tr-basis $X = \{X _{\Theta _{k}'}\colon \ \Theta
_{k} ^{\prime } \in {\rm Un}(\Theta _{k})\}$. Using notation as in
Proposition 4.1, put $\widetilde \Lambda _{\Theta _{k}'} =
\Lambda ((\Theta _{k} ^{\prime \prime }\Omega _{k})/\Omega _{k}; X
_{\Theta _{k}'})$, for each $\Theta _{k} ^{\prime } \in {\rm
Un}(\Theta _{k})$, where $\Theta _{k} ^{\prime \prime } = \Theta
_{k} ^{\prime } \cap (\Theta _{k,\chi '}\Theta _{k,{\rm Sol}})$.
Denote by $\widetilde \Theta _{k}$ the tensor compositum over
$\Omega _{k}$ of the fields $\widetilde \Lambda _{\Theta _{k}'}$,
$\Theta _{k} ^{\prime } \in {\rm Un}(\Theta _{k})$. It is not
difficult to see that $\Theta _{k}$ and $\Omega _{k}$ are
algebraically closed in $\widetilde \Theta _{k}$. At the same time,
one observes that $X _{\Theta _{k}'} \not\in N(O _{k}\Omega _{k}/F
_{k}\Omega _{k})$, for any finite extension $O _{k}$ of $\Theta
_{k}$ in $\Theta _{k,{\rm sep}}$ and any $F _{k} \in I(O
_{k}/\Theta _{k})$. Thus it follows that the conditions of Lemma 5.3
are fulfilled by $\Theta _{k}'\Omega _{k}/\Theta _{k}''\Omega _{k}$,
$X _{\widetilde \Theta _{k}'}$ and any $p \in \overline P$ dividing
$[\Theta _{k} ^{\prime }\colon \Theta _{k} ^{\prime \prime }]$.
Identifying $\Theta _{k,{\rm sep}}$ with its $\Theta
_{k}$-isomorphic copy in $\widetilde \Theta _{k,{\rm sep}}$, one
deduces from Proposition 4.1 and Lemmas 3.6, 5.3 and 2.2 that
Br$(\nabla _{k}\widetilde \Theta _{k}/\nabla _{k}) = \{0\}$, for
every finite extension $\nabla _{k}$ of $\Theta _{k}$ in $\Theta
_{k,{\rm sep}}$, and also, that $X _{\Theta _{k}'} \in N(\Theta _{k}
^{\prime \prime }\widetilde \Theta _{k}/\widetilde \Theta _{k})
\setminus N(\Theta _{k} ^{\prime }\widetilde \Theta _{k}/\Theta _{k}
^{\prime \prime }\widetilde \Theta _{k})$, for each $\Theta _{k}
^{\prime } \in {\rm Un}(\Theta _{k})$. Hence, by Lemma 2.2, $X
_{\Theta _{k}'} \not\in N(\Theta _{k} ^{\prime }\widetilde \Theta _{k}/\widetilde \Theta _{k})$, $\Theta _{k} ^{\prime } \in {\rm
Un}(\Theta _{k})$.
\par
Assuming that $\Delta _{k} = \widetilde \Theta _{k}$, if $\chi
^{\prime } = \chi $ or $\chi ^{\prime } \subseteq {\rm Sol}$, we
give the definition of $\Delta _{k}$ under the hypothesis that $\chi
\neq \chi ^{\prime } \neq {\rm Sol}$. Then $\chi ^{\prime }$ is
abelian closed and ${\rm Sol} \subset \chi ^{\prime }$. Let
$Z(\widetilde \Theta _{k}) = \widetilde \Theta _{k}\Sigma _{k,2}$, Un$(\widetilde \Theta _{k}) ^{\prime } = \{\widetilde \Theta _{k}
^{\prime } \in {\rm Gal}(\widetilde \Theta _{k})\colon \ \widetilde
\Theta _{k} ^{\prime } \subseteq Z(\widetilde \Theta _{k}),
G(\widetilde \Theta _{k} ^{\prime }/(\widetilde \Theta _{k} ^{\prime
} \cap \widetilde \Theta _{k,\chi })) \not\in {\rm Sol}\}$,
$\widetilde \Theta _{k} ^{\prime \prime } = \widetilde \Theta _{k}
^{\prime } \cap (\widetilde \Theta _{k,\chi }\widetilde \Theta
_{k,{\rm Sol}})$, for each $\widetilde \Theta _{k} ^{\prime } \in
{\rm Un}(\widetilde \Theta _{k}) ^{\prime }$, and let $\widetilde
\Omega _{k}/\widetilde \Theta _{k}$ be a rational extension with a
tr-basis $\widetilde X = \{\widetilde X _{\widetilde \Theta
_{k}'}\colon \ \widetilde \Theta _{k} ^{\prime } \in {\rm
Un}(\widetilde \Theta _{k}) ^{\prime }\}$. We define $\Delta _{k}$
to be the tensor compositum over $\widetilde \Omega _{k}$ of the
fields $\widetilde \Lambda _{\widetilde \Theta _{k}'} = \Lambda
((\widetilde \Theta _{k} ^{\prime \prime }\widetilde \Omega
_{k})/\widetilde \Omega _{k}; X _{\widetilde \Theta _{k}'})$,
$\widetilde \Theta _{k}' \in {\rm Un}(\widetilde \Theta _{k})
^{\prime })$. It is easily seen that $\widetilde \Theta _{k}$ and
$\widetilde \Omega _{k}$ are algebraically closed in $\Delta _{k}$.
Note also that $X _{\widetilde \Theta _{k}'} \not\in N(\widetilde O _{k}\widetilde \Omega _{k}/\widetilde F _{k}\widetilde \Omega _{k})$
in case $\widetilde O _{k}/\widetilde \Theta _{k}$ is a finite
extension and $\widetilde F _{k} \in I(\widetilde O _{k}/\widetilde
\Theta _{k})$. Hence, Lemma 5.3 applies to $\widetilde \Theta
_{k}'\widetilde \Omega _{k}/\widetilde \Theta _{k}''\widetilde
\Omega _{k}$, $X _{\widetilde \Theta _{k}'}$ and any $p \in
\overline P$ dividing $[\widetilde \Theta _{k} ^{\prime }\colon
\widetilde \Theta _{k} ^{\prime \prime }]$. Identifying $\widetilde
\Theta _{k,{\rm sep}}$ and $\widetilde \Omega _{k,{\rm sep}}$ with
their isomorphic copies in $\Delta _{k,{\rm sep}}$ (over $\widetilde
\Theta _{k}$ and $\widetilde \Omega _{k}$, respectively), one
deduces from Proposition 4.1 that Br$(\widetilde \nabla _{k}\Delta _{k}/\widetilde \nabla _{k}\widetilde \Omega _{k}) = \{0\}$ (and
Br$(\widetilde \nabla _{k}\Delta _{k}/\widetilde \nabla _{k}) =
\{0\}$), for each finite extension $\widetilde \nabla _{k}$ of
$\widetilde \Theta _{k}$ in $\widetilde \Theta _{k,{\rm sep}}$. Using
Galois theory and Lemmas 3.6, 5.3 and 2.2, one also proves the
following:
\par
\medskip
(5.3) (i) $X _{\Theta _{k}'} \in N(\Theta _{k} ^{\prime
\prime }\Delta _{k}/\Delta _{k}) \setminus N(\Theta _{k}
^{\prime }\Delta _{k}/\Theta _{k} ^{\prime \prime }\Delta
_{k}), \ \Theta _{k} ^{\prime } \in {\rm Un}(\Theta
_{k})$;
\par
(ii) $\widetilde X _{\widetilde \Theta _{k}'} \in N(\widetilde
\Theta _{k} ^{\prime \prime }\Delta _{k}/\Delta _{k}) \setminus
N(\widetilde \Theta _{k} ^{\prime }\Delta _{k}/\widetilde
\Theta _{k} ^{\prime \prime }\Delta _{k}), \ \widetilde
\Theta _{k} ^{\prime } \in {\rm Un}(\widetilde \Theta _{k}) ^{\prime
}$.
\par
\medskip
We are now in a position to finish the construction of $E _{k+1}$
and to show that the field $E = \cup _{n=1} ^{\infty } E _{n}$ is
quasilocal and has the properties required by Theorem 1.2 (i) and
(ii). Let $E _{k+1} = \Delta _{k}$, provided that $\chi \neq {\rm
Nil}$, and suppose further that $\chi = {\rm Nil}$. Denote by
$Z(\Delta _{k}) ^{\prime }$ the field $\Delta _{k,{\rm Sol}}$, if
$\chi ^{\prime } = {\rm Nil}$, and put $Z(\Delta _{k}) ^{\prime } =
(\Delta _{k}\Sigma _{k,2}) \cap \Delta _{k,{\rm Sol}}$, otherwise.
Fix a rational field extension $\widehat \Omega _{k}/\Delta _{k}$
with a tr-basis $\widehat X = \{\widehat X _{\Delta _{k}'}\}$,
indexed by all $\Delta _{k} ^{\prime } \in I(Z(\Delta _{k}) ^{\prime
}/\Delta _{k})$, for which $\Delta _{k} ^{\prime }/\Delta _{k}$
satisfies the conditions of Lemma 5.1. Taking as $E _{k+1}$ the
tensor compositum over $\widehat \Omega _{k}$ of the fields $\Lambda
((\Delta _{k} ^{\prime } \cap \Delta _{k,{\rm Nil}}).\widehat \Omega _{k}/\widehat \Omega _{k}; \widehat X _{\Delta _{k}'})$, $\widehat X
_{\Delta _{k}'} \in \widehat X$, we finish the construction of $E$.
It is easily verified that $E _{n+1}/E _{n}$ and $E/E _{n}$ are
regular, and that $E _{n}$, $\Lambda _{n}$, $W _{n}$, $\Theta _{n}$,
$\widetilde \Theta _{n}$ and $\Delta _{n}$ are algebraically closed
in $E _{n+1}$ and $E$, for all $n \in \hbox{\Bbb N}$. Hence, by
Galois theory and the observation preceding statement (5.1), $E$ has
the following property:
\par
\medskip
(5.4) For each $G \in {\rm Fin}$ and $j = 1,2 $, there exists $M
_{j} (G) \in {\rm Gal}(E) \cap I(\Sigma _{0,j}E/E)$, such that $G(M
_{j} (G)/E) \cong G$.
\par
\medskip
\noindent
Note also that $\rho _{E _{n}/E _{n+1}}$ maps $T _{n}$ bijectively
upon $T _{n+1}$, for every $n \in \hbox{\Bbb N}$, and each $L \in
{\rm Gal}(E)$ has a subfield $L _{k} \in {\rm Gal}(\Lambda _{k})$,
such that $G(L _{k}/\Lambda _{k}) \cong G(L/E)$, for some $k \in
\hbox{\Bbb N}$. When $L _{k} \subseteq Z _{k}$, this implies that $L
_{k}\Lambda _{n} \subseteq Z _{n}$ and Ker$({\rm Cor}_{(L
_{k}\Lambda _{n})/\Lambda _{n}}) \subseteq $ ${\rm Br}(L _{k}E _{n}/L
_{k}\Lambda _{n})$, for $n > k$. Therefore, $\rho
_{E _{0}/E}$ maps $T _{0} = T$ isomorphically on Br$(E)$, and
whenever $L \in {\rm Gal}(E)$ and $L \subseteq \cup _{n=1} ^{\infty
} Z _{n} \colon = Z _{\infty }$, Cor$_{L/E}$ is injective.
Similarly, it follows from Galois theory and the established
properties of $E _{n}$, $\Lambda _{n}$, $\Theta _{n}$ and $E _{n+1}$
that if $E ^{\prime }$ and $E _{p} ^{\prime }$ lie in Gal$(E)$ and
$E ^{\prime } \subseteq E _{p} ^{\prime } \subseteq E ^{\prime }
(p)$, for some $p \in \overline P$, then Br$(E _{p} ^{\prime }) _{p}
\cap {\rm Ker}({\rm Cor}_{E _{p}'/E'}) = \{0\}$. As Br$(E)$ is
divisible, these observations show that $E$ satisfies condition (c)
of Proposition 3.4, for $\Omega = E _{\rm sep}$ and each $p \in
\overline P$. Thus it turns out that $E$ is quasilocal and nonreal.
Applying now Corollary 3.5, one also concludes that $\rho _{E/R}$ is
surjective, for every finite extension $R$ of $E$ in $Z _{\infty }$.
Since, by (5.1), $E _{\chi } \subseteq Z _{\infty }$ in case $\chi
\neq {\rm Nil}$, the obtained results, combined with (1.1) (ii) and
Theorem 1.1 (i), prove Theorem 1.2 (i) and (ii).
\par
It remains to be seen that $E$ has the properties required by
Theorem 1.2 (iii) and (iv). Our argument goes along similar lines to
those drawn in the proof of Theorem 1.2 (i) and (ii), so we present
it omitting details. Let $M$ be a field lying in Gal$(E)$. Suppose
first that $G(M/(M \cap E _{\chi '})) \not\in {\rm Sol}$. Then $M =
M _{\kappa }E$, for some $\kappa \in \hbox{\Bbb N}$, $M _{k} \in
{\rm Un}(\Theta _{\kappa })$. Since $\Theta _{\kappa }$ is
algebraically closed in $E$, this implies that $G(M _{\kappa
}/\Theta _{\kappa }) \cong G(M/E)$. Let $L _{\kappa }$ be a minimal
subfield of $M _{\kappa }$ lying in Un$(\Theta _{k})$. Then it
follows from the definition of $E$, (5.3) (i) and (5.2) (i) that $X
_{L _{\kappa }}$ lies in $N(L \cap (E _{\chi '}E _{\rm Sol})/E)
\setminus N(L/L \cap (E _{\chi '}E _{\rm Sol}))$, where $L = L
_{\kappa }E$. Hence, by Lemma 2.2, $X _{L _{\kappa }} \not\in
N(L/E)$. As $N(M/E) \subseteq N(L/E)$ and, by Lemma 3.2, $N(L \cap
(E _{\chi '}E _{\rm Sol})/E) = N(L/E)N(M \cap (E _{\rm \chi '}E
_{\rm Sol})/E)$, this implies that $N(M/E) \neq N(M \cap (E _{\chi
'}E _{\rm Sol})/E)$. When $\chi ^{\prime } \neq {\rm Nil}$, the
obtained result proves Theorem 1.2 (iii), since then $E _{\rm Sol}
\subseteq E _{\chi '}$. Similarly, one gets from (5.3) (ii) and
(5.2) (ii) that $N(M/E) \neq N(M \cap (E _{\chi }E _{\rm Sol})/E)$
also in case $G(M/E) \in \chi ^{\prime }$, $G(M/(M \cap E _{\chi }))
\not\in {\rm Sol}$ and $M \subseteq \Sigma _{0,2}E$. Assume now that
$\chi = {\rm Nil}$, put $Z(E) ^{\prime } =$ $\cup _{n=1} Z(\Delta
_{n}) ^{\prime }$, and suppose that $M \subseteq Z(E) ^{\prime }$
and the pair $(M/E, \chi )$ is chosen as in Lemma 5.1. Then $G(M/E)
\in {\rm Sol}$ and our construction of $E$ guarantees that $N(M/E)
\neq N(M/E) _{\rm Nil}$. As each $M ^{\prime } \in {\rm Gal}(E)$,
with $G(M ^{\prime }/E) \not\in {\rm Nil}$ and $M ^{\prime } \subset
Z(E) ^{\prime }$, possesses a subfield satisfying the conditions of
Lemma 5.1 with respect to $E$, this enables one to deduce from (1.1)
(ii) and Theorem 1.1 (ii) that $N(M ^{\prime }/E) \neq N(M ^{\prime
}/E) _{\rm Ab}$. Thus Theorem 1.2 (iii) is proved. For the proof of
Theorem 1.2 (iv), drop the assumption that $\chi = {\rm Nil}$, let
$\chi ^{\prime } \neq \chi $ and take fields $M _{1} (G)$ and $M
_{2} (G)$ as in (5.4), for an arbitrary $G \in \chi ^{\prime }
\setminus \chi $. As shown in the process of proving Theorem 1.2
(i)-(ii), then $\Sigma _{0,1}E \cap E _{\chi '} \subseteq Z
_{\infty }$ and $N(M _{1} (G)/E) = N(M _{1} (G)/E) _{\rm Ab}$. On
the contrary, Theorem 1.2 (ii), the proof of Theorem 1.2 (iii) and
the inclusion ${\rm Ab} \subset \chi $ indicate that $N(M _{2}
(G)/E) \neq N(M _{2} (G)/E) _{\chi } = N(M _{2} (G)/E) _{\rm Ab}$,
so Theorem 1.2 is proved.
\par
\medskip
{\bf Remark 5.4.} With assumptions being as in Theorem 1.2, suppose
that $E _{0}$ is infinite of cardinality $d \ge {\rm max}\{d
_{p}\colon \ p \in \overline P\}$, where $d _{p}$ is the dimension
of the subgroup $\{t _{p} \in T\colon \ pt _{p} = 0\} \subseteq T$
as a vector space over $\hbox{\Bbb F} _{p}$, for each $p \in
\overline P$. Analyzing the proofs of (4.3) (iv) and Theorem 1.2,
one concludes that our construction of the quasilocal field $E$ can
be specified so that $d$ equals the cardinality of $E$ as well as
the ranks of $G(E (p)/E)$ and of any Sylow pro-$p$-subgroup $G _{p}$
of $G _{E}$ as pro-$p$-groups, for each $p \in \overline P$.
Applying [7, II, Lemma 3.3], one also sees that whenever $E/E _{0}$,
$T$ and $\chi = {\rm Fin}$ are related as in Theorem 1.2, $d _{p}$
equals the dimension of the (continuous) cohomology group $H ^{2} (G
_{p}, \hbox{\Bbb F} _{p})$ as an $\hbox{\Bbb F} _{p}$-vector space.
Moreover, by the same lemma, then the dimension of $H ^{2} (G(E
(p)/E), \hbox{\Bbb F} _{p})$ is equal to $d _{p}$ or zero, depending
on whether or not $E _{0}$ contains a primitive $p$-th root of
unity.
\par
\vskip0.6truecm
{\centerline{\bf 6. Topological interpretation of Theorem 1.2}
\smallskip
\centerline{\bf and applications to Brauer groups of fields}}
\par
\medskip
Let $E$ be a field, char$(E) = q \ge 0$, $D(E)$ and $R(E)$ the
divisible and the reduced parts of $E ^{\ast }$, respectively. For
each (nonempty group) formation $\Psi \subseteq {\rm Fin}$, put $N
_{\Psi } (E) = \cap _{M \in {\rm Gal}(E)} N(M/E) _{\Psi }$. It is
easily obtained that $E ^{\ast }/N _{\Psi } (E)$ can be canonically
viewed as a topological group (totally disconnected, see [28,
Theorem 9]). Specifically, $N _{\Psi } (E) = D(E)$, if $E$ is
perfect and quasilocal, Br$(E) _{p} \neq \{0\}$, for every $p \in
(\Pi (E) \setminus \{q\})$, and ${\rm Met} \cap {\rm Fin}_{\Pi (E)}
\subseteq \Psi $, where ${\rm Fin}_{\Pi (E)}$ is the class of all $G
\in {\rm Fin}$ of orders not divisible by any $p \in \overline P
\setminus \Pi (E)$ (see [7, II, Sect. 2]). With the same
assumptions on $E$, the equality $N _{\Psi } (E) = D(E)$ holds also
when $P(E) \setminus \{q\} = \Pi (E) \setminus \{q\}$ and ${\rm Ab}
\cap {\rm Fin}_{\Pi (E)} \subseteq \Psi $. Hence, by the isomorphism
$R(E) \cong E ^{\ast }/D(E)$, $\Psi $ induces on $R(E)$ a structure
$T _{\Psi }$ of a topological group. Let $\Psi _{1}$ be another
formation of this kind. It follows from Lemma 3.2 that the
topologies $T _{\Psi }$ and $T _{\Psi _{1}}$ are equivalent ($T
_{\Psi } \sim T _{\Psi _{1}}$) if and only if $N(M/E) = N(M/E)
_{(\Psi \cap \Psi _{1})}$ when $M \in {\rm Gal}(E)$ and $M \subseteq
(E _{\Psi } \cup E _{\Psi _{1}})$. Clearly, this occurs if and only
if $T _{\Psi } \sim T _{(\Psi \cap \Psi _{1})}$ and $T _{\Psi _{1}}
\sim T _{(\Psi \cap \Psi _{1})}$, so solving the equivalence problem
for $T _{\Psi }$ and $T _{\Psi _{1}}$ reduces to the case in which
$\Psi \subset \Psi _{1}$. Consider now the set $\Omega _{\rm ac}$ of
all pairs $(\chi , \chi ^{\prime })$ of distinct subclasses of ${\rm
Fin}$ satisfying the conditions of Theorem 1.2. A binary sequence
$\bar a = a _{i,j}\colon \ (i, j) \in \Omega _{\rm ac}$, is called
multiplicative, if $a _{k,l}a _{l,m} = a _{k,m}$ for every $(k, l)$
and $(l, m)$ lying in $\Omega _{\rm ac}$. When $P(E) \setminus \{q\}
= \Pi (E) \setminus \{q\}$, such a sequence $\bar a(E)$ can be
canonically attached to $E$, putting $a _{i,j} (E) = 1$, if $T _{i}
\sim T _{j}$, and $a _{i,j} (E) = 0$, otherwise. By Theorem 1.1 (i),
$a _{i,j} (E) = 1$, for all $(i, j) \in \Omega _{\rm ac}$, provided
that $\rho _{E/M}$ is surjective, for each $M \in {\rm Gal}(E)$.
This applies to the presently known, and conjecturally, to all
perfect fields with LCFT in the sense of Neukirch-Perlis (see [8,
Proposition 3.3 and Remark 3.4 (ii)]). Conversely, when $\bar a$ is
multiplicative and $T$ is a divisible abelian torsion group, one
obtains by modifying the proof of Theorem 1.2 that there exists a
quasilocal perfect field $F$, such that:
\par
\medskip
(6.1) (i) Br$(F) \cong T$ and all $G \in {\rm Fin}$ are realizable
as Galois groups over $F$; in particular, $\Pi (F) = P(F) =
\overline P$.
\par
(ii) $\bar a = \bar a (F)$ and $\rho _{F _{i}/M _{(i)}}$ is
surjective when $a _{i,j} = 1$ and $M _{(i)} \in {\rm Gal}(F _{i})
\cap I(F _{j}/E)$.
\par
\medskip \noindent
The groups $D(E)$ and $R(E)$ are related with group formations as
above, also if $E$ is SQL and almost perfect (in the sense of [7]).
\par
\medskip
{\bf Remark 6.1.} Let $R _{\rm sim}$ be a system of representatives
of the isomorphism classes of finite simple groups, $\Sigma _{\rm
sim} = \{S \subseteq R _{\rm sim}\colon \ {\rm Ab} \cap R _{\rm
sim} \subseteq S\}$, and for each $S \in \Sigma _{\rm sim}$, let
$\widetilde S$ be the class of those $G \in {\rm Fin}$, whose simple
quotients are isomorphic to groups from $S$. It follows from the
Jordan-H\"older theorem that the correspondence $S \to \widetilde
S\colon \ S \in \Sigma _{\rm sim}$, is injective and maps $\Sigma
_{\rm sim}$ into the set of abelian closed classes.
\par
\medskip
Now we turn our attention to the Brauer groups of the basic types of
PQL-fields. First we present two results which substantially
generalize Theorem 1.2 (i) and (ii). They also complement Theorem
1.1 (i), [8, Corollary 6.2] and an observation made by M. Auslander
(see [31, Ch. II, Sect. 3.1]) about the class of fields with trivial
Brauer groups.
\par
\medskip
{\bf Proposition 6.2.} {\it Let $E _{0}$ be a field, $L _{0}$ a
Galois extension of $E _{0}$ in $E _{0,{\rm sep}}$, $S$ a nonempty
set of profinite groups, and let $T$ and $T _{0}$ satisfy the
conditions of Theorem 1.2. Then there exists a regular field
extension $E/E _{0}$ with the following properties:}
\par
(i) $E$ {\it is quasilocal, {\rm Br}$(E) \cong T$, $\rho _{E
_{0}/E}$ maps $T _{0}$ injectively into {\rm Br}$(E)$, and all $G
\in S$ are realizable as Galois groups over $E$;}
\par
(ii) {\it A finite extension $R$ of $E$ in $E _{\rm sep}$ lies in
$I(L _{0}E/E)$ if and only if $\rho _{E/R}$ is surjective; when $R
\not\in I(L _{0}E/E)$ and $p \in \overline P$ does not divide
$[\widetilde R\colon (R \cap L _{0}E)]$, where $\widetilde R
\subseteq E _{\rm sep}$ is the normal closure of $R$ over $E$, {\rm
Br}$(R) _{p}$ properly includes the image of {\rm Br}$(E) _{p}$
under $\rho _{E/R}$.}
\par
\medskip
{\it Proof.} The existence of $E$ is proved constructively, and in
this respect, our proof is very similar to the one of Theorem 1.2
(i) and (ii), in the special case where $\chi = {\rm Fin}$. In the
first place, it becomes clear that one may additionally assume that
$T _{0} = T$, all $G \in S$ are realizable as Galois groups over $E
_{0}$, and each $H \in {\rm Ab}$ has an isomorphic copy $s _{H} \in
S$. In this setting, $E$ is obtained as a union $\cup _{n=0}
^{\infty } \Lambda _{n} = \cup _{n=0} ^{\infty } W _{n} = \cup
_{n=0} ^{\infty } \Theta _{n} = \cup _{n=0} ^{\infty } E _{n}$ of
inductively defined field towers, such that $E _{n} \subseteq
\Lambda _{n} \subseteq W _{n} \subseteq \Theta _{n} \subseteq E
_{n+1}$, for every index $n$. The extensions $\Lambda _{n}/E _{n}$
and $\Theta _{n}/W _{n}$ are defined in exactly the same way as in
the proof of Theorem 1.2 (whence they are regular), $E _{n+1}/\Theta
_{n}$ are rational of tr-degree $2$, and $W _{n}$ is $\widetilde
\otimes _{\Lambda _{n}}$ of the fields $W(M _{n}/\Lambda _{n})$,
where $M _{n}$ runs across Gal$(\Lambda _{n}) \cap I(L _{0}\Lambda
_{n}/\Lambda _{n})$. This ensures that $E$ has the properties
required by Proposition 6.2 (i) as well as the surjectivity of $\rho
_{E/R}$, for every finite extension $R$ of $E$ in $L _{0}E$. In
particular, the regularity of $E/E _{0}$ and the additional
conditions satisfied by $E _{0}$ and $S$ enable one to deduce from
Galois theory that $C(\Phi (\pi )/\Phi )$ contains infinitely many
elements of order $\pi $, for each $\Phi \in I(E/E _{0})$ and $\pi
\in \overline P$. It remains for us to prove the latter assertion of
Proposition 6.2 (ii), so we assume further that $E$, $R$ and
$\widetilde R$ satisfy its conditions. Fix some $p \in \overline P$
not dividing $[\widetilde R\colon V]$, where $V = R \cap L _{0}E$, put 
$U _{n} = U \cap \Theta _{n,{\rm sep}}$, for each index $n$ and $U \in
I(\widetilde R/E)$, and define Im$_{p} (V/R)$ and Ker$_{p} (R/V)$ in
accordance with the proof of Corollary 3.5. Using the RC-formula and
Proposition 2.3 (ii), one obtains that it suffices to establish the
inequality Im$_{p} (V/R) \neq {\rm Br}(R) _{p}$ under the extra
hypothesis that $I(R/E) \setminus \{R\} \subseteq I(L _{0}E/E)$. It
follows from Proposition 2.3 (ii) and the choice of $p$ that Br$(R)
_{p} \cong {\rm Im}_{p} (V/R) \oplus {\rm Ker}_{p} (R/V)$ (see also
[27, Sect. 13.4]), so our assertion is equivalent to the one that
Ker$_{p} (R/V) \neq \{0\}$. Clearly, $\widetilde R = \widetilde R
_{k}E$, for some integer $k \ge 0$ and $\widetilde R _{k} \in {\rm
Gal}(\Theta _{k})$. As $E/\Theta _{k}$ is regular, $G(\widetilde R
_{k}/\Theta _{k})$ and $G(\widetilde R/E)$ are canonically
isomorphic, and for convenience, they will be further identified.
Fix a tr-basis $X _{k}, Y _{k}$ of $E _{k+1}/\Theta _{k}$ and a
primitive element $\varphi _{k}$ of $\widetilde R _{k}/\Theta _{k}$.
Observe that $E _{\rm sep}$ contains as a subfield a cyclic
extension $\Theta _{k} ^{\prime }$ of $\Theta _{k} (Y _{k})$ of
degree $p$, such that $\Theta _{k} ^{\prime } \cap \widetilde R _{k}
(Y _{k}) = \Theta _{k} (Y _{k})$. This can be easily deduced from
Galois theory and the fact that the maximal subgroup of $C(\Theta
_{k} (Y _{k}) (p)/\Theta _{k} (Y _{k}))$ of exponent $p$ is
infinite. Regarding now the groups $\widetilde R _{k} (X
_{k})/\widetilde R _{k} (X _{k}) ^{\ast p}$ and $_{p} {\rm
Br}(\widetilde R _{k}E _{k+1})$ as modules over the group algebra
$\hbox{\Bbb F} _{p} [G(\widetilde R/E)]$, and considering the images
of the element $X _{k} - \varphi _{k}$ under the action of
$G(\widetilde R/E)$, one obtains without difficulty that $\hbox{\Bbb
F} _{p} [G(\widetilde R/E)]$ is $\hbox{\Bbb F} _{p}$-isomorphic to
submodules of $\widetilde R _{k} (X _{k})/\widetilde R _{k} (X _{k})
^{\ast p}$ and $_{p} {\rm Br}(\widetilde R _{k}E _{k+1})$. Let
$[R\colon V] = r$ and $\{\xi _{j}\colon \ j = 1, \dots , r\}$ be a
system of representatives of the right co-sets of $G(\widetilde
R/R)$ in $G(\widetilde R/V)$. The embeddability of $\hbox{\Bbb F}
_{p} [G(\widetilde R/E)]$ in $_{p} {\rm Br}(\widetilde R _{k}E
_{k+1})$ indicates that $_{p} {\rm Br}(\widetilde R _{k}E _{k+1})$
contains an element $\delta $, such that $[\sum _{j=1} ^{r} (\xi
_{j} - 1)]\delta $ lies in Ker$(\widetilde R _{k}E _{k+1}/V _{k}E
_{k+1}) \setminus {\rm Ker}(\widetilde R _{k}E _{k+1}/R _{k}E
_{k+1})$. This implies that Ker$_{p} (R _{k}E _{k+1}/V _{k}E
_{k+1}) \neq \{0\}$. At the same time, it follows from Galois theory
and the regularity of $E/\Theta _{k}$ that $U _{n} = U _{k}\Theta
_{n}$, $V _{n} = R \cap L _{0}\Theta _{n}$ and $G(\widetilde R
_{n}/\Theta _{n}) \cong G(\widetilde R/E)$, for $n \ge k$ and $U \in
I(\widetilde R/E)$. Hence, by Proposition 4.6 (i), the extra
hypothesis on $R$ and the noted properties of the construction of
$E$, Br$(R _{n}) _{p} \cap {\rm Br}(R _{n+1}/R _{n}) \subseteq {\rm
Im }_{p} (V _{n}/R _{n})$, for every index $n \ge k$. In view of the
RC-formula and the choice of $p$, these observations show that $\rho
_{R _{k}/R}$ maps Ker$_{p} (R _{k}/V _{k})$ injectively into
Ker$_{p} (R/V)$. In particular, it turns out that Ker$_{p} (R/V)
\neq \{0\}$, so Proposition 6.2 is proved.
\par
\medskip
Our next result, applied to the formation $\Psi $ of supersolvable
groups $W \in {\rm Fin}$, and to a set $S$ of profinite groups
containing an isomorphic copy of the symmetric group Sym$_{4}$,
proves the existence of a quasilocal field $E$, such that Sym$_{4}$
is realizable as a Galois group over $E$, and each $M \in {\rm
Gal}(E)$ with $G(M/E) \cong {\rm Sym}_{4}$ satisfies both the
condition of Corollary 3.10 and the one that $\rho _{E/\Phi }$ is
not surjective, for any $\Phi \in {\rm Gal}(E)$ including $M$. This
indicates that the assertions of Theorem 1.2 (i), (ii) and (iv)
cannot simultaneously be true, if $\chi $ and $\chi ^{\prime }$ are
replaced by $\Psi $ and ${\rm Sol}$, respectively.
\par
\medskip
{\bf Proposition 6.3.} {\it Let $T$ and $S$ satisfy the conditions of
Proposition 6.2, and let $\Psi \subseteq {\rm Fin}$ be a
formation. Then there exists a quasilocal field $E$, such that {\rm
Br}$(E) \cong T$, every $G \in S$ is realizable as a Galois group
over $E$, and a finite extension $R$ of $E$ in $E _{\rm sep}$ lies
in $I(E _{\Psi }/E)$ if and only if $\rho _{E/R}$ is surjective.}
\par
\medskip
{\it Proof.} The field $E$ can be obtained as an extension of an
arbitrary fixed field $E _{0}$, which is a union $E = \cup _{n=0}
^{\infty } E _{n}$ of an inductively defined tower, such that $E
_{n}/E _{n-1}$ has the properties required by Proposition 6.2 with
respect to $E _{n-1,\Psi }$, for each $n \in \hbox{\Bbb N}$.
\par
\medskip
Proposition 2.3 (ii) and the following result, combined with [15,
Theorem 23.1], describe the isomorphism classes of Brauer groups of
PQL-fields.
\par
\medskip
{\bf Proposition 6.4.} {\it Let $E _{0}$ be a formally real field,
${\rm Odd} \subset {\rm Fin}$ the class of groups of odd orders,
$\Psi \subseteq {\rm Odd}$ a formation, and $T$ an abelian torsion
group, such that the $p$-components $T _{p}\colon \ p \in \overline
P \setminus \{2\}$, are divisible, and $T _{2}$ of order $2$.
Suppose also that $T _{0}$ is a subgroup of {\rm Br}$(E _{0})$
embeddable in $T$ and with $T _{0,2} = \{0\}$. Then there exists a
field extension $E/E _{0}$ with the following properties:}
\par
(i) $E$ {\it is formally real and {\rm PQL}, $E _{0} (2)$ contains
as a subfield the algebraic closure of $E _{0}$ in $E$, {\rm Br}$(E)
\cong T$ and $T _{0} \cap {\rm Br}(E/E _{0}) = \{0\}$;}
\par
(ii) {\it Finite extensions of $E$ in $E _{\rm Odd}$ are
$p$-quasilocal, for every $p \in \overline P \setminus \{2\}$, and
all $G \in {\rm Odd}$ are realizable as Galois groups over $E$;}
\par
(iii) {\it A finite extension $R$ of $E$ in $E _{\rm Odd}$ is
included in $E _{\Psi }$ if and only if {\rm Br}$(R) _{p} \subseteq
{\rm Im}(E/R)$, for every $p \in \overline P \setminus \{2\}$; if
$R \not\in I(E _{\Psi }/E)$, then {\rm Br}$(R) _{p'} \cap {\rm
Im}(E/R) \neq {\rm Br}(R) _{p'}$, for any $p ^{\prime } \in
\overline P \setminus \{2\}$, $p ^{\prime } \not\vert [R\colon E 
_{\Psi } \cap R]$.}
\par
\medskip
{\it Proof.} We obtain $E$ as a union $\cup _{n=1} ^{\infty } E
_{n}$ of a field tower defined inductively as in the proof of
Proposition 6.3. Omitting the details, note that $E$ is constructed by
considering ${\rm Odd}$, $\overline P \setminus \{2\}$ and the sums
of the subgroups $T _{0,p}\colon \ p > 2$, Br$(E _{n}) _{p}\colon \
p > 2$, instead of ${\rm Fin}$, $\overline P$, $T _{0}$ and Br$(E
_{n})$, respectively, extending at the end of the inductive step the
obtained field by a real closure in its maximal $2$-extension. The
correctness of the construction is proved essentially as in the
proof of Proposition 6.2; the Artin-Schreier theory ensures that $E
_{n}$, $n \in \hbox{\Bbb N}$, and $E$ are formally real, since the
function field of each of the Brauer-Severi varieties $V$ used for
constructing $E$ embeds over the field of definition $\Phi (V)$ of
$V$ into a field $R(V)$ that is rational over some extension $\Phi
^{\prime } (V)$ of $\Phi (V)$ of odd degree.
\par
\medskip
{\bf Remark 6.5.} In the setting of Proposition 6.4, when $\Psi =
\{1\}$, the proof of Theorem 1.2 enables one to modify the
construction of $E$ so as to satisfy the inequality $N(M/E) \neq
N(M/E) _{\rm Ab}$, for each $M \in {\rm Gal}(E)$ with $G(M/E) \in
{\rm Odd} \setminus {\rm Nil}$.
\par
\medskip
Our next result complements [8, Theorems 1.2, 1.3 and Proposition 3.6]
as follows:
\par
\medskip
{\bf Corollary 6.6.} {\it Let $T$ satisfy the condition of some of 
Propositions 6.3 or 6.4, $\Psi = {\rm Fin}$ in the former case, $\Psi
= {\rm Odd}$ in the latter one, and let ${\rm Sup}(T) = \{p \in
\overline P\colon \ T _{p} \neq \{0\}\}$. Then there is a strictly
{\rm PQL}-field $F$, such that {\rm Br}$(F) \cong T$ and every $G
\in \Psi $, for which the index $\vert G\colon [G, G]\vert $ has no
divisor $p \in \overline P \setminus {\rm Sup}(T)$, is realizable as
a Galois group over $F$. When $T$ is divisible and ${\rm Sup}(T) =
\overline P$, $F$ can be chosen from the class of {\rm SQL}-fields.}
\par
\medskip
{\it Proof.} Our latter conclusion follows at once from Theorem
1.2, and in case ${\rm Sup}(T) = \overline P$, the former one is
contained in Theorem 1.2 and Proposition 6.4. Assume now that ${\rm
Sup}(T) \neq \overline P$, put $\overline {\rm Sup}(T) = \overline P
\setminus {\rm Sup}(T)$, consider a field $F _{0}$ with Br$(F _{0})
\cong T$, and let $T _{0} = \{0\} \subset {\rm Br}(F _{0})$. Then
one can take as $F$ the union $\cup _{n=0} F _{n} ^{\prime } = \cup
_{n=0} F _{n} ^{\prime \prime } = \cup _{n=1} F _{n}$ of fields
defined inductively so as to satisfy the following conditions, for
each index $n \ge 0$:
\par
\medskip
(6.2) (i) Br$(F _{n} ^{\prime }) \cong T$ and $F _{n} ^{\prime }/F
_{n}$ has the properties required by Theorem 1.2 or Proposition 6.4.
More precisely, $\rho _{F _{0}/F _{n}'}$ is an isomorphism and, for
every $p \in \overline P \setminus \{2\}$ and each finite
extension $R _{n}$ of $F _{n} ^{\prime }$ in $F _{n,\Psi } ^{\prime
}$, $\rho _{F _{n}'/R _{n}}$ maps Br$(F _{n} ^{\prime }) _{p}$
surjectively on Br$(R _{n}) _{p}$.
\par
(ii) $F _{n} ^{\prime \prime }$ is the compositum of the fields $F
_{n} ^{\prime } (p)\colon \ p \in \overline {\rm Sup}(T)$.
\par
(iii) $F _{n+1} = F _{n} ^{\prime \prime }$, if $T$ is divisible,
and $F _{n+1}$ is a real closure of $F _{n} ^{\prime \prime }$ in $F
_{n} ^{\prime \prime } (2)$, otherwise.
\par
\medskip \noindent
When $T$ and $F _{0}$ have the properties required by Theorem 1.2
(i)-(ii), for $\chi = {\rm Fin}$, one may use only a countable
iteration of (6.2) (ii), with $(F _{n}, F _{n+1})$ instead of $(F
_{n} ^{\prime }, F _{n} ^{\prime \prime })$, for any $n$.
\par
\medskip
{\bf Corollary 6.7.} {\it Let $E$ be a field and $T$ an abelian
torsion group with the properties required by Theorem 1.2 (i),
Proposition 6.4 or Corollary 6.6. Assume $E$ has a Henselian 
valuation $v$. Then the value group $v(E)$ of $(E, v)$ is divisible 
and every $D \in d(E)$ is defectless with respect to $v$.}
\par
\medskip
{\it Proof.} Consider first a Henselian valued quasilocal field $(F, w)$, 
such that $w(F) = qw(F)$, Br$(F) _{q} \neq \{0\}$, char$(\widehat F) = 
q > 0$, and $q \in \Pi (\widehat F)$, where $\widehat F$ is the residue 
field of $(F, w)$. By [20, Theorem 3.16], $F$ and $\widehat F$ are 
nonreal fields. For each finite extension $L/F$, denote by $w(L)$ the 
value group, and by $\widehat L$ the residue field of $L$ with respect 
to its unique (up-to an equivalence) valuation $w _{L}$ extending $w$. 
It is well-known that $w(F)$ is a subgroup of $w(L)$ of index $\vert 
w(L)\colon w(F)\vert \le [L\colon F]$ (cf. [21, Ch. XII, Proposition 12]), 
and that $w(L)$ is isomorphic to a totally ordered subgroup of $w(F)$. 
The equality $w(F) = qw(F)$ implies that $q$ does not divide $\vert 
w(L)\colon w(F)\vert $, and yields $w(L) = qw(L)$ as well. We show that 
$L/F$ is defectless with respect to $w$ (and $w _{L}$), provided that $L 
\subset F _{\rm sep}$. The inequality cd$_{q} (G _{\widehat F}) > 0$ 
and [19, Theorem 2.8 (a)] guarantee the existence of an inertial finite 
extension $\Phi /F$, such that $q \not\vert [\Phi \colon F]$ and $q \in 
P(\widehat \Phi )$. At the same time, it is easily seen that Br$(\Phi /F) 
\cap {\rm Br}(F) _{q} = \{0\}$ (cf. [27, Sect. 13.4]), whence the nontriviality 
of Br$(F) _{q}$ is preserved by Br$(\Phi ) _{q}$. Moreover, it follows from 
Ostrowski's theorem (cf. [12]) that our assertion holds if and only if finite 
extensions of $\Phi $ in $\Phi _{\rm sep}$ are defectless with respect to 
$v _{\Phi }$, so the preceding observations reduce its proof to the 
special case in which $q \in P(\widehat F)$. Then the quasilocal property 
of $F$ and the nontriviality of Br$(F) _{q}$ imply that $N(L _{1}/F) \neq F 
^{\ast }$, for every inertial cyclic extension $L _{1}$ of $F$ in $F (q)$. 
Hence, by the lifting property of $w$, applied to a norm form of $L _{1}/F$ 
with coefficients in the valuation ring of $(F, w)$, the assumption on 
$w(F)$ ensures that $N(\widehat L _{1}/\widehat F) \neq 
\widehat F ^{\ast }$. In view of [27, Sect. 15.1, Proposition b] and [19, 
Theorem 2.8 (a)], this proves the existence of an inertial $F$-algebra 
$\Delta _{q} \in d(F)$ of index $q$. Observing that every extension of 
$F$ embeddable in $\Delta _{q}$ is inertial, one obtains from the 
quasilocal property of $F$ that there are no immediate cyclic extensions 
of $F$ of degree $q$. Since, by Witt's theorem, Br$(F) _{q}$ is divisible, 
and by [42, Theorem 2] and Galois theory, $q \in P(\widehat R)$, for 
every finite extension $R/F$, this observation makes it easy to deduce 
the claimed defectlessness of $L/F$ (in case $L \subset F _{\rm sep}$) 
from Ostrowski's theorem and well-known formulae about valuation 
prolongations. As shown in [36], the obtained result implies that every 
$\Delta \in d(F)$ is defectless. When $w(F)$ is divisible, this means that 
$\Delta $ is inertial over $F$. Also, it follows from the Ostrowski-Draxl 
theorem [12], that our conclusions remain valid, if $(F, v)$ is a 
Henselian valued field, such that char$(\widehat F) = 0$. Thus the 
latter conclusion of Corollary 6.7 turns out to be a consequence of the 
former one.
\par
Our objective now is to establish the divisibility of $v(E)$. In the first 
two cases, this follows directly from the former assertion of [8,  (2.3)]. In 
the third one, one obtains from Galois theory that if $E$ is nonreal, then 
simple groups $\Sigma \in {\rm Fin} \setminus {\rm Ab}$ are realizable 
as Galois groups over $E$, so our assertion again is implied by the 
noted part of [8, (2.3)]. Assume further that $E$ is formally real, Br$(E) 
\cong T$, ${\rm Sup}(T) \neq \overline P$ and $E$ has the properties 
required by Corollary 6.6. By [20, Theorem 3.16], then the residue field 
$\widehat E$ of $(E, v)$ is formally real, whence char$(\widehat E) = 0$. 
In addition, it is well-known that if $G \in {\rm Fin}$, $p \in \overline P$ 
and $\nu $ is the order of $G$, then $G$ embeds in Aut$(P(G))$, where 
$P(G) \in {\rm Ab}$ is of exponent $p$ and order $p ^{\nu }$. When $G 
\not\in {\rm Ab}$ and $p \not\vert \nu $, this implies that $G$ has an 
irreducible representation over $\hbox{\Bbb F} _{p}$ of dimension $\ge 
2$. Hence, by Galois theory and the assumptions on $E$, for each $p 
\in \overline P \setminus \{2\}$, there exists $M _{p} \in {\rm Gal}(E) \cap 
I(E _{\rm Odd}/E)$, such that $p \not\vert [M _{p}\colon E]$ and $G(M 
_{p} (p)/M _{p})$ is of rank $\ge 2$ as a pro-$p$-group (one may put 
$M _{p} = E$, if $p \in {\rm Sup}(T)$). Let now $v _{p}$ be the 
prolongation of $v$ on $M _{p}$. By the proof of [8, (2.3)], then $v _{p} 
(M _{p}) = pv _{p} (M _{p}) = 2v _{p} (M _{p})$ ($M _{p}$ is formally 
real), and since $v(E)$ is a subgroup of $v _{p} (M _{p})$ of index $\le 
[M _{p}\colon E]$, this means that $v(E) = pv(E)$, $p \in \overline P$, 
which completes our proof.
\par
\medskip
In conclusion, we use Theorem 1.2 (i) for describing, up-to an
isomorphism, the abelian torsion groups that can be realized as
reduced parts of Brauer groups of absolutely stable fields $K$
possessing equicharacteristic Henselian valuations $v$, such that
$v(K)$ are totally indivisible (i.e. with $v(K) \neq pv(K)$, for
every $p \in \overline P$). This enables one to construct various
new examples of Brauer groups that are not simply presentable (see
[16, Lemma 1 and pages 492-493]). Before giving the description,
note that the residue field $\widehat K$ is quasilocal, and that a
Henselian discrete valued field $(L, w)$ is absolutely stable,
provided that $\widehat L$ is quasilocal and perfect (see [7, I,
Proposition 2.3 and the beginning of Sect. 8]).
\par
\medskip
{\bf Proposition 6.8.} {\it An abelian torsion group $\Theta $ is
isomorphic to a maximal reduced subgroup of {\rm Br}$(K)$, for an
absolutely stable field $K = K(\Theta )$ with a Henselian valuation
$v$ such that {\rm char}$(K) = {\rm char}(\widehat K)$ and $v(K)$ is
totally indivisible if and only if the $p$-component $\Theta _{p}$
decomposes into a direct sum of cyclic groups of the same order $p
^{n _{p}}$, for each $p \in \overline P$.}
\par
\medskip
{\it Proof.} The necessity of the conditions on $\Theta $ has been
proved in [7, II, Sect. 3], so we show here only their sufficiency.
Let $\widetilde \Theta $ be a divisible hull of $\Theta $, $\Pi
(\Theta ) = \{p \in \overline P\colon \ n _{p} > 0\}$, $\overline
{\hbox{\Bbb Q}}$ an algebraic closure of $\hbox{\Bbb Q}$, and
$\varepsilon _{n} \in \overline {\hbox{\Bbb Q}}$ a primitive $n$-th
root of unity, for each $n \in \hbox{\Bbb N}$. Denote by $F _{0}$
the extension of $\hbox{\Bbb Q}$ generated by the set $\{\varepsilon
_{p ^{n _{p}}}\colon \ p \in \Pi (\Theta )\} \cup \{\varepsilon _{p
^{\nu }}, \nu \in \hbox{\Bbb N}\colon \ p \in \overline P \setminus
\Pi (\Theta )\}$. Take an extension $E/F _{0}$ in accordance with
Theorem 1.2 so that Br$(E) \cong \widetilde \Theta $, and consider a
Henselian discrete valued field $(K(\Theta ), v)$ with a residue
field $F _{0}$-isomorphic to $E$. By Scharlau's generalization of
Witt's decomposition theorem [31], then Br$(K(\Theta ))$ $\cong {\rm
Br}(E) \oplus C _{E}$ (see also [38, (3.10)]), so it follows from
(2.5) that the reduced part of Br$(K(\Theta ))$ is isomorphic to
$\Theta $. Thus Proposition 6.8 is proved.
\par
\medskip
Note finally that the groups $\Theta $ singled out by Proposition
6.8 are those realizable as reduced parts of Brauer groups of
Henselian discrete valued absolutely stable fields. The sufficiency
of the conditions on $\Theta $ is shown by the proof of Proposition
6.8, and their necessity follows from [7, II, Lemma 3.2], [38,
(3.10)] and the divisibility of Br$(F) _{q}$ and $C(F (q)/F)$, for
any field $F$ of characteristic $q \neq 0$ (Witt, see [11, Sect.
15]). On the other hand, each sequence $\{p ^{h _{p}}\colon \ h _{p}
\in (\hbox{\Bbb N} \cup \{\infty \}), \ p \in \overline P\}$ with $h
_{2} \neq \infty $ and $p ^{h _{p}} \colon = \infty $, $h _{p} =
\infty $, equals the sequence $\{e _{p} (K)\colon \ p \in \overline
P\}$ of exponents of reduced parts of Br$(K) _{p}$, for some
Henselian discrete valued stable field $(K, v)$ (see the reference
at the end of [10, Sect. 3]).
\par
\vskip1.6truecm
\centerline{\bf References}
\par
\medskip
[1] A.A. Albert, Modern Higher Algebra. Univ. of Chicago
Press, XIV, Chicago, Ill., 1937.
\par
[2] E. Artin, J. Tate, Class Field Theory. Benjamin, New
York-Amsterdam, 1968.
\par
[3] A. Blanchet, Function fields of generalized Brauer-Severi
varieties. Comm. Algebra 19 (1991), No 1, 97-118.
\par
[4] L. Br\" ocker, Characterization of fans and hereditarily
Pythagorean fields. Math. Z. 151 (1976), 149-173.
\par
[5] I.D. Chipchakov, On nilpotent Galois groups and the scope
of the norm limitation theorem in one-dimensional abstract local
class field theory. In: Proc. of ICTAMI 05, Alba Iulia, Romania,
15.9-18.9, 2005; Acta Univ. Apulensis 10 (2005), 149-167.
\par
[6] I.D. Chipchakov, Class field theory for strictly
quasilocal fields with Henselian discrete valuations. Manuscr.
math. 119 (2006), 383-394.
\par
[7] I.D. Chipchakov, On the residue fields of Henselian
valued stable fields. I, J. Algebra 319 (2008), 16-49; II, C.R.
Acad. Bulg. Sci 60 (2007), No 5, 471-478.
\par
[8] I.D. Chipchakov, One-dimensional abstract local class
field theory. Preprint, arXiv:math/0506515v5 [math.RA].
\par
[9] I.D. Chipchakov, Norm groups and class fields of formally
real quasi-local fields. Preprint, arXiv:math/0508019v1 [math.RA].
\par
[10] I.D. Chipchakov, Embeddings in quasilocal fields and
computations in Brauer groups of arbitrary fields. C.R. Acad. 
Bulg. Sci. 61 (2008), 1229-1238.
\par
[11] P.K. Draxl, Skew Fields. Lond. Math. Soc. Lecture Notes,
Cambridge etc., Cambridge Univ. Press, 1983.
\par
[12] P.K. Draxl, Ostrowski's theorem for Henselian valued skew 
fields. J. Reine Angew Math. 354 (1984), 213-218. 
\par
[13] B. Fein, D.J. Saltman, M. Schacher, Embedding problems for
finite dimensional division algebras. J. Algebra 167 (1994), No 3,
588-626.
\par
[14] B. Fein, M. Schacher, Relative Brauer groups, I. J.
Reine Angew. Math. 321 (1981), 179-194.
\par
[15] L. Fuchs, Infinite Abelian Groups. V. I and II, Pure and
Applied Mathematics, Academic Press, IX and XI, New York-London,
1970.
\par
[16] A.W. Hales, Abelian groups as Brauer and character
groups. Abelian group theory, Proc. 3rd Conf., Oberwolfach/FRG 1985,
487-507 (1987).
\par
[17] K. Iwasawa, Local Class Field Theory. Iwanami Shoten,
Japan, 1980 (Japanese: Russian transl. in Mir, Moscow, 1983).
\par
[18] G. Karpilovsky, Topics in Field Theory. North-Holland
Math. Studies, 155, North Holland, Amsterdam etc., 1989.
\par
[19] B. Jacob and A. Wadsworth, Division algebras over Henselian 
fields. J. Algebra 128 (1990), 126-179. 
\par
[20] T.Y. Lam, Orderings, valuations and quadratic forms. Conf.
Board Math. Sci. Regional Conf. Ser. Math. No 52, Amer. Math. Soc.,
Providence, RI, 1983.
\par
[21] S. Lang, Algebra. Addison-Wesley Publ. Comp., Mass., 1965.
\par
[22] P. Mammone, A. Merkurjev, On the corestriction of the $p
^{n}$-symbol. Isr. J. Math. 76 (1991), 73-80.
\par
[23] A.S. Merkurjev, Index reduction formula. J. Ramanujan Math.
Soc. 12 (1997), 49-95.
\par
[24] A.S. Merkurjev, I.A. Panin, A.R. Wadsworth, Index reduction
formulas for twisted flag varieties. I, K-Theory 10 (1996), 517-596;
II, K-Theory 14 (1998), 101-196.
\par
[25] A.S. Merkurjev, A.A. Suslin, $K$-cohomology of Brauer-Severi
varieties and the norm residue homomorphism. Izv. Akad. Nauk SSSR 46
(1982), 1011-1046 (Russian: Engl. transl. in Math. USSR Izv. 21
(1983), 307-340).
\par
[26] J. Neukirch, R. Perlis, Fields with local class field theory.
J. Algebra 42 (1976), 531-536.
\par
[27] R. Pierce, Associative Algebras. Graduate Texts in Math. 88,
Springer-Verlag, New York-Heidelberg-Berlin, 1982.
\par
[28] L.S. Pontrjagin, Topological Groups, 4-th Ed. Nauka, Moscow,
1984 (Russian).
\par
[29] L. Redei, Algebra, v. 1. Akademiai Kiado, Budapest, 1967.
\par
[30] D.J. Saltman, The Schur index and Moody's theorem. K-Theory 7
(163), No 3, 309-332.
\par
[31] W. Scharlau, \"Uber die Brauer-Gruppe eines
Hensel-K\"orpers. Abh. Math. Semin. Univ. Hamb. 33 (1969),
243-249.
\par
[32] J.-P. Serre, Cohomologie Galoisienne. Lecture Notes in Math. 5,
Springer-Verlag, Berlin-Heidelberg-New York, 1965.
\par
[33] J.-P. Serre, Local Fields. Graduate Texts in Mathematics 88,
Springer-Verlag, Berlin-Heidelberg-New York, 1979.
\par
[34] L.A. Shemetkov, Formations of Finite Groups. Nauka, Moscow,
1978 (in Russian).
\par
[35] J.-P. Tignol, On the corestriction of central simple algebras.
Math. Z. 194 (1987), No 2, 267-274.
\par
[36] I.L. Tomchin and V.I. Yanchevskij, On defects of valued division 
algebras. Algebra i Analiz 3 (1991), No 3, 147-164 (Russian: Engl. 
transl. in St. Petersbg. Math. J. 3 (1992), No 3, 631-647).
\par
[37] M. Van den Bergh, A. Schofield, Division algebra coproducts of
index $n$. Trans. Amer. Math. Soc. 341 (1994), 505-518.
\par
[38] A.R. Wadsworth, Valuation theory on finite dimensional division
algebras. In: F.-V. Kuhlmann (ed.) et al., Valuation Theory and its
applications, I. Proc. of international conference and workshop,
Univ. of Saskatchewan, Saskatoon, Canada, 28.7-11.8, 1999.
\par
[39] W.C. Waterhouse, Profinite groups are Galois groups. Proc.
Amer. Math. Soc. 42 (1974), 639-640.
\par
[40] E. Weiss, Cohomology of Groups. Pure and Applied Mathematics
34, Academic Press, New York-London, 1969.
\par
[41] G. Whaples, Generalized local class field theory. II. Existence
Theorem. Duke Math. J. 21 (1954), 247-255.
\par
[42] G. Whaples, Algebrtaic extensions of arbitrary fields. Duke Math. J. 
24 (1957), 201-204.
\par
\bye